\documentclass[VANCOUVER,STIX1COL]{WileyNJD-v2}

\articletype{RESEARCH ARTICLE}%

\received{}
\revised{}
\accepted{}

\usepackage{lipsum}
\usepackage{bm}
\usepackage{amsmath}
\usepackage{graphicx}
\usepackage{epstopdf}

\ifpdf
  \DeclareGraphicsExtensions{.eps,.pdf,.png,.jpg}
\else
  \DeclareGraphicsExtensions{.eps}
\fi

\numberwithin{theorem}{section}

\usepackage{amsopn}
\DeclareMathOperator{\diag}{diag}

\usepackage{euscript}
\usepackage{relsize}
\usepackage{mathdots}
\usepackage{graphicx}
\usepackage{indentfirst}
\usepackage{empheq}
\usepackage{multirow}
\usepackage{thmtools}
\usepackage{thm-restate}

\usepackage{tikz}

\DeclareMathOperator*{\argmax}{arg\,max}
\DeclareMathOperator*{\argmin}{arg\,min}

\newcommand{\row}[1]{#1 ,\,\boldsymbol{:}\,}
\newcommand{\vecrow}[1]{#1}
\newcommand{\col}[1]{\,\boldsymbol{:}, #1}

\def\rank{\mathop{\mathrm{rank}}}

\def\colspace{\mathop{\mathrm{colspace}}}

\def\Sigminus{\bfW}
\newcommand\winverse[2]{\left(#1 \right)^{\dag}_{#2}}
\newcommand\inverse[1]{#1^{\dag}}
\def\i0{\tau}
\def\Si0{\dot \bfs}
\def\Ai0{\dot \bfa}
\def\Proj{\mathbf{\Pi}}

\def\fullop{H_{\tau}}

\newcommand{\Input}{\hspace*{\algorithmicindent} \textbf{Input}:\ }

\usepackage{euscript}

\newcommand{\spC}{\mathbb{C}}

\newcommand{\spR}{\mathbb{R}}

\newcommand{\spT}{\mathbb{T}}

\newcommand{\tsA}{\mathsf{A}}
\newcommand{\tsB}{\mathsf{B}}

\newcommand{\tsR}{\mathsf{R}}
\newcommand{\tsS}{\mathsf{S}}

\newcommand{\tsX}{\mathsf{X}}
\newcommand{\tsY}{\mathsf{Y}}
\newcommand{\tsZ}{\mathsf{Z}}

\newcommand{\bfA}{\mathbf{A}}
\newcommand{\bfB}{\mathbf{B}}
\newcommand{\bfC}{\mathbf{C}}
\newcommand{\bfG}{\mathbf{G}}
\newcommand{\bfL}{\mathbf{L}}
\newcommand{\bfX}{\mathbf{X}}
\newcommand{\bfF}{\mathbf{F}}
\newcommand{\bfY}{\mathbf{Y}}
\newcommand{\bfZ}{\mathbf{Z}}
\newcommand{\bfM}{\mathbf{M}}

\newcommand{\bfQ}{\mathbf{Q}}
\newcommand{\bfT}{\mathbf{T}}

\newcommand{\bfW}{\mathbf{W}}
\newcommand{\bfI}{\mathbf{I}}
\newcommand{\bfR}{\mathbf{R}}
\newcommand{\bfS}{\mathbf{S}}
\newcommand{\bfU}{\mathbf{U}}
\newcommand{\bfO}{\mathbf{O}}
\newcommand{\bfzero}{\mathbf{0}}

\newcommand{\bfJ}{\mathbf{J}}

\newcommand{\rmT}{\mathrm{T}}

\newcommand{\rmF}{\mathrm{F}}
\newcommand{\rmW}{\mathrm{W}}
\newcommand{\rmR}{\mathrm{R}}

\newcommand{\calD}{\mathcal{D}}
\newcommand{\calC}{\mathcal{C}}
\newcommand{\calT}{T}
\newcommand{\calF}{\mathcal{F}}

\newcommand{\calL}{\mathcal{L}}

\newcommand{\calG}{\mathcal{G}}
\newcommand{\calK}{\mathcal{K}}
\newcommand{\calI}{\mathcal{I}}
\newcommand{\calQ}{\mathcal{Q}}

\newcommand{\calW}{\mathcal{W}}
\newcommand{\calZ}{\mathcal{Z}}

\newcommand{\bfa}{\mathbf{a}}
\newcommand{\bfb}{\mathbf{b}}
\newcommand{\bfc}{\mathbf{c}}
\newcommand{\bfd}{\mathbf{d}}
\newcommand{\bfe}{\mathbf{e}}

\newcommand{\bfp}{\mathbf{p}}
\newcommand{\bfq}{\mathbf{q}}

\newcommand{\bfs}{\mathbf{s}}

\newcommand{\bfu}{\mathbf{u}}
\newcommand{\bfv}{\mathbf{v}}

\newcommand{\bfx}{\mathbf{x}}
\newcommand{\bfy}{\mathbf{y}}

\newcommand{\bt}{\begin{theorem}}
\newcommand{\et}{\end{theorem}}
\newcommand{\bl}{\begin{lemma}}
\newcommand{\el}{\end{lemma}}
\newcommand{\bp}{\begin{proposition}}
\newcommand{\ep}{\end{proposition}}
\newcommand{\bc}{\begin{corollary}}
\newcommand{\ec}{\end{corollary}}

\newcommand{\bd}{\begin{definition}\rm}
\newcommand{\ed}{\end{definition}}
\newcommand{\bex}{\begin{example}\rm}
\newcommand{\eex}{\end{example}}
\newcommand{\br}{\begin{remark}\rm}
\newcommand{\er}{\end{remark}}

\newcommand{\btbh}{\begin{table}[!ht]}
\newcommand{\etb}{\end{table}}
\newcommand{\bfgh}{\begin{figure}[!ht]}
\newcommand{\efg}{\end{figure}}

\newcommand{\bea}{\begin{eqnarray*}}
\newcommand{\eea}{\end{eqnarray*}}
\newcommand{\be}{\begin{eqnarray}}
\newcommand{\ee}{\end{eqnarray}}

\def\rank{\mathop{\mathrm{rank}}}

\makeatletter
\def\adots{\mathinner{\mkern2mu\raise\p@\hbox{.}
\mkern2mu\raise4\p@\hbox{.}\mkern1mu
\raise7\p@\vbox{\kern7\p@\hbox{.}}\mkern1mu}}
\newcommand{\l@abcd}[2]{\hbox to\textwidth{#1\dotfill #2}}
\makeatother

\def\unit{\mathfrak{i}}

\def\code#1{\texttt{#1}}

\begin{document}

\title{Fast and stable modification of the Gauss-Newton method for low-rank signal estimation}

\author[1]{Nikita Zvonarev}

\author[1]{Nina Golyandina}

\address[1]{\orgdiv{Faculty of Mathematics and Mechanics}, \orgname{St.Petersburg State University}, \orgaddress{\state{Universitetskaya nab. 7/9, St.Petersburg}, \country{Russia}}}

\corres{*Nina Golyandina, SPbSU, Universitetskaya nab. 7/9, St.Petersburg, 199034, Russia. \email{n.golyandina@spbu.ru}}

\abstract[Summary]{The weighted nonlinear least-squares problem for low-rank signal estimation is considered. The problem of constructing a numerical solution that is stable and fast for long time series is addressed. A modified weighted Gauss-Newton method, which can be implemented through the direct variable projection onto a space of low-rank signals, is proposed.
For a weight matrix which provides the maximum likelihood estimator of the signal in the presence of autoregressive noise of order $p$ the computational cost of iterations is $O(N r^2 + N p^2 + r N \log N)$ as $N$ tends to infinity, where $N$ is the time-series length, $r$ is the rank of the approximating time series. Moreover, the proposed method can be applied to data with missing values, without increasing the computational cost. The method is compared with state-of-the-art methods based on  the  variable projection approach in terms of floating-point numerical stability and computational cost.
}


\maketitle

\section{Introduction}\label{sec:intro}
In this study we consider  the `signal plus noise' observation scheme:
$$
x_n=s_n+\epsilon_n, \;\;\; n=1,2, \ldots, N.
$$
Denote by
$\tsX = (x_1, \ldots, x_N)^\rmT$ , $\tsS = (s_1, \ldots, s_N)^\rmT$  and $\bm\epsilon = (\epsilon_1, \ldots, \epsilon_N)^\rmT $
the vectors of observations,  signal values and errors respectively.
We will refer to vectors of observations in $\spR^N$ as time series (or shortly series, since the observations are not necessarily
temporal; e.g., they can be spatial).

We assume that the signal $\tsS$ can  be written in the parametric form as a finite sum
    \begin{equation}
    \label{eq:model}
    s_n = \sum_{k=1}^d P_{m_k}(n)\exp(\alpha_k n) \sin(2\pi \omega_k n + \phi_k),
    \end{equation}
 where   $P_{m_k}(n)$ are polynomials in $n$ of degree $m_k$.
 In signal processing applications, the signal in the form  \eqref{eq:model} is usually a sum of sine waves \cite{Cadzow1988}
 or a sum of damped sinusoids \cite{Markovsky2008}.
  The problem of estimating the unknown signal values $s_n$ is as important as the problem of estimating the parameters in the explicit form \eqref{eq:model}. We will be concentrated on the signal estimation.

  Certainly, there are many numerical methods for estimating signals in the form \eqref{eq:model} in the `signal plus noise' observation scheme. Our aim is to construct a method, which is accurate and fast for large series lengths $N$ and also is competitive. It would be a great advantage if the constructed method could deal with time series with gaps keeping algorithms' efficiency characteristics. The case of autoregressive noise is of special interest.

     We will solve the problem of signal estimation in a class of signals, which is wider than that given in \eqref{eq:model} and has a different parameterization from the explicit one.
To describe the class of signals, we need some definitions.
       The rank of a signal $\tsS$  is defined as follows.
    For a given integer $L$, $1<L<N$, called the window length, we define the embedding operator $\calT_{L}:\; \spR^{N} \to \spR^{L \times (N-L+1)}$,
    which maps $\tsS$ into a Hankel matrix, by
    	\begin{equation}
\label{eq:embedding}
	\calT_{L}(\tsS) = \begin{pmatrix}
	s_1 & s_2 & \hdots & s_{N-L+1} \\
	s_2 & s_3 & \hdots & \vdots \\
	\vdots & \vdots & \hdots & s_{N-1} \\
	s_{L} & s_{L+1} & \hdots & s_{N}
	\end{pmatrix}.
	\end{equation}
The columns of $\calT_{L}(\tsS)$ are sequential lagged vectors; this is why $\calT_{L}(\tsS)$ is often called the \emph{trajectory matrix} of $\tsS$.
		We say that the signal $\tsS$ has rank $r < N/2$ if $\rank \calT_{r+1}(\tsS) = r$.
 It is known that $\rank \calT_{r+1}(\tsS) = r$ if and only if $\rank \calT_{L}(\tsS) = r$ for any $L$ such that $\min(L, N-L+1) > r$ (see \cite[Corollary 5.1]{Heinig1984} for the proof).
	
For a sufficiently large time-series length $N$, the signal in the form \eqref{eq:model} has rank $r$, which is determined by the parameters $m_k$, $\alpha_k$ and
$\omega_k$ \cite{Zvonarev.Golyandina2021}.
For example, the signal with values $s_n$ has rank $r=2$ for a sum of two exponentials $s_n=c_1 \exp(\alpha_1 n) + c_2 \exp(\alpha_2 n)$, a sine wave $s_n=c \sin(2 \pi \omega n + \phi)$, where $0<\omega<0.5$, or a linear function $s_n =a n + b$. Thus, the model fixes the rank $r$ but does not fix the form of the signal,
    i.e., the number of terms and the orders of the polynomials in \eqref{eq:model}.

   The considered model of signals, where the Hankel matrix $\calT_L(\tsS)$ is rank-deficient, is one of the standard models in many areas, signal processing \cite{Cadzow1988,Tufts1993}, speech recognition \cite{Dendrinos1991}, control theory and linear systems \cite{Markovsky2008,Markovsky2019} among others.

\medskip
    Let us state the optimization problem for estimating the signal.
    Denote $\calD_r$ the set of series of rank $r$.
	Since the set $\calD_r$ is not closed, we will seek for the solution in its closure $\overline{\calD_r}$. It is well-known that $\overline{\calD_r}$ consists of series of rank not larger than $r$ (this result can be found in \cite[Remark 1.46]{iarrobino1999power} for the complex case; the real-valued case is considered in \cite{Zvonarev.Golyandina2021}).

Thus, in what follows, we study the weighted least-squares (WLS) problem
	\begin{equation}
    \label{eq:wls}
	\tsY^\star = \argmin_{\tsY \in \overline{\calD_r}} \| \tsX - \tsY \|_{\bfW},
	\end{equation}
where $\bfW$ is a symmetric positive-definite weight matrix, $\|\tsZ\|_{\bfW}^2 = \tsZ^\rmT \bfW \tsZ$.
If noise $\bm\epsilon$ is Gaussian with covariance matrix $\bm\Sigma$ and zero mean, the WLS estimate with the weight matrix $\bfW = \bm\Sigma^{-1}$ is the maximum likelihood estimator (MLE). The same is true if $\bfW$ is scaled by a constant. In the case when noise is an autoregressive process of order $p$ (AR($p$)), the matrix $\bfW$ is $(2p+1)$-diagonal \cite[p.534]{Shaman1975}.

Although the common case is the case of a positive definite matrix $\bfW$, the problem \eqref{eq:wls}, where
 $\bfW$ is positive semi-definite, is of considerable interest. For example, the case of a diagonal matrix $\bfW$
 with several zero diagonal elements corresponds to the problem of low-rank approximation for time series with missing values if the noise is white. Let us consider the case of a general weight matrix and time series with missing values.  Let a symmetric positive definite matrix $\bfW_0$ be given for the whole time series including gaps. Then the weight matrix $\bfW$ is constructed from $\bfW_0$ by setting the columns and rows with indices equal to entries of missing values to zero.
 Note that if $\bfW$ is not positive-definite, then $\|\cdot\|_{\bfW}$ is semi-norm and the problem \eqref{eq:wls} may become ill-posed. In particular, the topology of $\overline{\calD_r}$ is not consistent with the semi-norm and
 therefore the minimum in \eqref{eq:wls} should be changed to infimum, which can be not achieved at time series in $\overline{\calD_r}$.
 	
\textbf{Different approaches for solving  \eqref{eq:wls}.}
The optimization problem \eqref{eq:wls} is  non-convex  with many local minima \cite{Ottaviani2014lraexact}.
The problem \eqref{eq:wls} is commonly considered as a structured (more precisely, Hankel) low-rank approximation problem (SLRA, HSLRA) \cite{Chu2003,Markovsky2006,Markovsky2019}.
A well-known subspace-based method for solving \eqref{eq:wls}
is called `Cadzow iterations' \cite{Cadzow1988} and belongs to the class
of alternating-projection methods.
The method of Cadzow iterations can be extended to a class of oblique Cadzow iterations in the norm,
which differs from the Euclidean norm \cite{Gillard2016}.
The Cadzow method has two drawbacks: first, the properties of the limiting point of the Cadzow iterations are unknown \cite{Andersson2013}
and second, it tries to solve the problem \eqref{eq:wls} with a weight matrix which generally differs
from the given $\Sigminus$.
Therefore, it is not optimal (the method does not provide the MLE), even for the case of white Gaussian noise \cite{DeMoor1994}.
The reason is that the problems are commonly stated in SLRA as matrix approximation problems, while the original problem
\eqref{eq:wls} is stated in terms of time series.

Many methods have been proposed to solve  HSLRA, including the Riemannian SVD \cite{DeMoor1994}, Structured total least-norm \cite{Lemmerling2000}, Newton-like iterations \cite{Schost2014}, proximal iterations \cite{Condat.Hirabayashi2015}, symbolic computations \cite{Ottaviani2014lraexact}, stochastic optimization \cite{Gillard2013}, fixed point iterations \cite{Andersson2013}, a penalization approach \cite{Ishteva.etal2014}.

Since we consider the problem of WLS time series approximation, which generally differs from the problem of matrix approximation due to different weights (see e.g. \cite{Zvonarev.Golyandina2017}), let us use as a benchmark the effective and general approach of Markovsky and Usevich \cite{Usevich2012,Usevich2014}, which is based on the variable projection principle \cite{Golub.Pereyr2003} combined with the Gauss-Newton method for solving the optimization subproblem; we call it VPGN. The method of \cite{Usevich2012,Usevich2014} is able to deal with the problem in the form \eqref{eq:wls}, i.e., exactly with the given weight matrix; moreover, it is elaborated in general form for a wide class of structured matrices and at the same time its iteration complexity scales linearly with the length of data for a class of weight matrices. Thus, the VPGN method can be considered as a start-of-art method of low-rank time series approximation.
Nevertheless, the approach has a couple of disadvantages. First, the Cholesky factorization is used for solving least-squares subproblems to obtain a fast algorithm; unfortunately, this squares the condition number (more stable  decompositions like QR factorization are slower). Then, the VPGN method is efficient only if the inverse of the weight matrix is banded.
Note that the approach of Markovsky and Usevich can be applied to the case of rank-deficient matrices $\bfW$ in \cite{Markovsky2013missing} and \cite[Section 4.4]{Markovsky2019}. However, it is not clear how to implement the algorithm, which proposed in these papers, effectively in terms of computational costs.
In \cite{Zvonarev.Golyandina2021}, a version S-VPGN, which improves the stability of several steps of VPGN is considered; however, it does not overcome these drawbacks.

\textbf{The proposed approach.} Let us consider another approach to solving the problem \eqref{eq:wls}; the considered approach is similar to VPGN but different.
It is discussed in Section~\ref{sec:lrr} that each time series $\tsS$ in $\overline\calD_r$ is characterized by a vector $\bfa\in \spR^{r+1}$, which provides the coefficients of a generalized linear recurrence relation (GLRR) governing the time series, i.e. $\bfa^\rmT \calT_{r+1}(\tsS)$ is the zero vector. For each $\bfa$, we can consider the space $\calZ(\bfa)$ of signals governed by the GLRR with the given coefficients. Iterative algorithms that use the variable projection method for solving the problem \eqref{eq:wls} include the projection onto $\calZ(\bfa)$ as a subproblem.

In this paper, we propose to overcome the drawbacks of the methods of \cite{Usevich2012,Usevich2014} and \cite{Zvonarev.Golyandina2021} in the following manner. First, the proposed modification helps to avoid computing the pseudoinverse of the Jacobian matrix (compare \eqref{eq:gauss_simple} and \eqref{eq:iterGNfinal}).
Then, as well as in \cite{Zvonarev.Golyandina2021} (and unlike \cite{Usevich2014}), the projection is calculated directly onto the space $\calZ(\bfa)$ and is not obtained through projecting on its orthogonal complement $\calQ(\bfa)$.
Finally, for calculating the projection, we use fast algorithms with improved numerical stability (the compensated Horner scheme, see Remark~\ref{rem:hornerscheme}).
As a result, the proposed method sometimes can be just slightly slower and is much faster in many real-life scenarios, but also is more stable (see Section~\ref{sec:comparison} with the comparison results).
Moreover, the algorithm with direct projections onto the space $\calZ(\bfa)$ can be extended to the case of a degenerate weight matrix $\bfW$
 (in particular, to the case of missing values) without loss of effectiveness (see Section~\ref{sec:missing}); compare with that in \cite{Markovsky2013missing}, where projections to the subspace $\calQ(\bfa)$ are used and the computational cost considerably increases for degenerate weight matrices.

\textbf{Structure of the paper}.
In  Section~\ref{sec:parameterization} we briefly discuss the parameterization of $\calD_r$ and its properties and introduce necessarily notation.
In Section~\ref{sec:optim} we describe the known (VPGN) and the new proposed (MGN) iterative methods
for solving the optimization problem~\eqref{eq:wls}. The algorithm VPGN is described in the way different from
that in \cite{Usevich2014}, since the description in \cite{Usevich2014} is performed for general SLRA problems, whereas we consider a particular case of time series.
Section~\ref{sec:MGNand VPGN} presents the algorithms with the implementations of VPGN and MGN.
In Section~\ref{sec:comparison} we compare computational costs and numerical stability of the VPGN and MGN algorithms.
Section~\ref{sec:conclusion} concludes the paper.

\textbf{Main notation}.
In this paper, we use lowercase letters ($a$,$b$,\ldots) and also $L$, $K$, $M$, $N$ for scalars, bold lowercase letters ($\bfa$,$\bfb$,\ldots) for vectors, bold uppercase letters ($\bfA$,$\bfB$,\ldots) for matrices, and
the calligraphic font for sets. Formally, time series are vectors; however, we use the uppercase sans serif font ($\tsA$,$\tsB$,\ldots) for time series to distinguish them from ordinary vectors.
Additionally, $\bfI_{M}\in \spR^{M\times M}$ is the identity matrix, $\bm{0}_{M \times k}$ denotes the $M \times k$ zero matrix, $\bm{0}_M$ denotes the zero vector in $\spR^M$,
$\bfe_i$ is the $i$-th standard basis vector.

Denote $\bfb_{\vecrow{\calC}}$ the vector consisting of
the elements of a vector $\bfb$ with indices in a set $\calC$.
For matrices, denote $\bfB_{\row{\calC}}$ the matrix consisting
of rows of a matrix $\bfB$ with indices in $\calC$ and $\bfB_{\col{\calC}}$
the matrix consisting
of columns of a matrix $\bfB$ with indices in $\calC$.

Finally, we put a brief list of main common symbols and acronyms.\\
LRR is a linear recurrence relation.\\
GLRR($\bfa$) is a generalized LRR with the coefficients given by $\bfa$.\\
$\calD_r$ is the set of time series of rank $r$.\\
$\overline{\calD_r}$ is the set of time series of rank not larger than $r$.\\
$\calZ(\bfa) \in \spR^N$ is the set of time series of length $N$ governed by the minimal GLRR($\bfa$);
$\bfZ(\bfa)$ is a matrix consisting of its basis vectors.\\
$\calQ(\bfa)$ is the orthogonal complement to $\calZ(\bfa)$; $\bfQ(\bfa)$ is the matrix consisting of its special basis vectors (see Section~\ref{subsec:subspace_approach}).\\
$\bfW\in \spR^{N\times N}$ is a weight matrix.\\
$\winverse{\bfF}{\bfW}$ is the weighted pseudoinverse matrix; $\inverse{\bfF}$ stands for $\winverse{\bfF}{\bfI_N}$.\\
$\bfJ_{S}$ is the Jacobian matrix of a map $S$.\\
$\calT_M: \spR^N \rightarrow \spR^{M\times (N-M+1)}$ is the embedding operator, which constructs the $M$-trajectory matrix.\\
$\fullop$: $\spR^{r} \to \spR^{r+1}$ is the operator, which inserts $-1$ in the position $\tau$.\\
$\Proj_{\calL,\bfW}$ is the $\bfW$-orthogonal projection onto $\calL$, $\Proj_{\bfL,\bfW}$ is the $\bfW$-orthogonal projection onto the column space $\colspace(\bfL)$; if $\bfW$ is the identity matrix, it is omitted in the notation.\\
$S_{\tau}^\star(\Ai0) = \Proj_{\calZ(\fullop(\Ai0)), \bfW}\tsX$, where $\Ai0\in \spR^r$.

	\section{Parameterization of low-rank series}
    \label{sec:parameterization}
    In this section, we introduce the notations and describe the results of \cite{Zvonarev.Golyandina2021} which we will use further.

    \subsection{Generalized linear recurrence relations}
    \label{sec:lrr}

    It is well known \cite[Theorem 3.1.1]{Hall1998} that for sufficiently large $N$,
    a time series in the form \eqref{eq:model} satisfies a linear recurrence relation (LRR) of some order $m$:
    \begin{equation}
    \label{eq:lrr}
    s_n = \sum_{k=1}^m b_k s_{n-k}, n = m+1, \ldots, N; b_m\neq 0.
    \end{equation}
    One time series can be governed by many different LRRs. The LRR of minimal order $r$ (it is unique) is called minimal.
    If $r<N/2$, the corresponding time series has rank $r$. The minimal LRR uniquely defines the form of \eqref{eq:model} and the parameters $m_k$, $\alpha_k$, $\omega_k$.

    The relations \eqref{eq:lrr} can be expressed in vector form as
    $\bfa^\rmT \calT_{m+1}(\tsS) = \bfzero_{N-m}^\rmT$, where $\bfa = (b_m, \ldots, b_1, -1)^\rmT \in \spR^{m+1}$.
    The vector $\bfa$ corresponding to the minimal LRR ($m=r+1$) and the first $r$  values of the series $\tsS$ uniquely determine the whole series  $\tsS$.
    Therefore, $r$ coefficients of an LRR of order $r$ and $r$ initial values
    ($2r$ parameters altogether) can be chosen as parameters of a series of rank $r$.
    However, this parameterization does not describe the whole set. For example, $\tsS=(1,1,1,1,1,2)^\rmT$ has rank 2 and does not satisfy an LRR of order $m<N-1$.

Let us generalize LRRs.
We say that \emph{a time series satisfies a generalized LRR (GLRR) of order $m$} if $\bfa^\rmT \calT_{m+1}(\tsS) = \bfzero_{N-m}^\rmT$ for some non-zero $\bfa \in \spR^{m+1}$;
we call this linear relation GLRR($\bfa$). If $a_{m+1}=-1$, the same relation is called LRR($\bfa$).
As well as for LRRs, the minimal GLRR can be introduced.
The difference between a GLRR  and an ordinary LRR is that the last coefficient in the GLRR is not necessarily non-zero and
therefore the GLRR does not necessarily set a recurrence.
However, at least one of the coefficients of the GLRR should be non-zero.

Let us demonstrate the difference between LRRs and GLRRs by an example. Let $\tsS = (s_1,\ldots,s_N)^\rmT$ be a signal and $\bfa = (a_1, a_2, a_3)^\rmT$.
Then GLRR($\bfa$) and LRR($\bfa$) mean the same:
$a_1 s_i + a_2 s_{i+1} + a_3 s_{i+2} = 0$ for $i=1,\ldots,N-2$. For LRR($\bfa$), we state that $a_3  = -1$ (or just not equal to 0). Then this linear relation becomes a recurrence relation since $s_{i+2} = a_1 s_i + a_2 s_{i+1}$.
For GLRR($\bfa$), we assume that some of $a_i$ is not zero (or equal to $-1$). It may be $a_1$ or $a_2$ or $a_3$. For example, $\tsS=(1,1,1,1,1,2)^\rmT$ satisfies the GLRR($\bfa$) with $\bfa = (1,-1,0)^\rmT$.

Thus, we consider the parameterization with the help of GLRR($\bfa$); in fact, the same approach is used in \cite{Usevich2012, Usevich2014}. In what follows, we assume that $2r<N$.

The following properties clarify the difference between the spaces $\calD_r$ and its closure $\overline{\calD_r}$:
(a)
$\overline{\calD_r} = \{\tsY: \exists \bfa\in \spR^{r+1}, \bfa\neq \bfzero_{r+1}: \bfa^\rmT \calT_{r+1}(\tsS) = \bfzero_{N-r}^\rmT\}$ or, equivalently,
$\tsY \in \overline{\calD_r}$ if and only if there exists a GLRR($\bfa$) of order $r$, which governs $\tsY$;
(b)
$\tsY \in \calD_r$ if and only if there exists a GLRR($\bfa$) of order $r$, which governs $\tsY$, and this GLRR is minimal.

\subsection{Subspace approach}
\label{subsec:subspace_approach}
	Let  $\calZ(\bfa)$, $\bfa\in \spR^{r+1}$, be the space of time series of length $N$ governed by the GLRR($\bfa$);
that is, $\calZ(\bfa) = \{\tsS: \bfa^\rmT \calT_{r+1}(\tsS) = \bfzero_{N-r}^\rmT \}$.
Therefore $\overline{\calD_r} = \bigcup \limits_{\bfa} \calZ(\bfa)$.

Define the operator $\bfQ:\,\spR^{r+1} \to \spR^{M \times (N - r)}$ as
\begin{equation}\label{op:Q}
\big(\bfQ(\bfa)\big)^\mathrm{T} = \begin{pmatrix}
a_1 & a_2 & \dots & \dots & a_{r+1} & 0 & \dots & 0 \\
0 & a_1 & a_2 & \dots &  \dots & a_{r+1} & \ddots & \vdots \\
\vdots & \ddots  & \ddots & \ddots & \ddots & \ddots & \ddots & 0 \\
0 & \dots & 0 & a_1 & a_2 & \ddots & \ddots & a_{r+1} \\
\end{pmatrix}.
	\end{equation}
Then the other convenient form of $\calZ(\bfa)$ is $\calZ(\bfa) = \{\tsS: \bfQ^\rmT(\bfa) \tsS = \bfzero_{N-r} \}$.

The following notation will be used below: $\calQ(\bfa) = \colspace(\bfQ(\bfa))$ and denote $\bfZ(\bfa)$ a matrix whose column vectors form a basis of $\calZ(\bfa)$.

\subsection{Parameterization}
\label{sec:param}
Let us describe the details of the considered parameterization.
Consider a series $\tsS_0\in \calD_r$, which satisfies a minimal GLRR($\bfa_0$) of order $r$ defined by a non-zero vector $\bfa_0 = (a_1^{(0)}, \ldots, $ $a_{r+1}^{(0)})^\rmT$.
Let us fix $\i0$ such that $a^{(0)}_{\i0} \neq 0$. Since GLRR($\bfa_0$) is invariant to multiplication by a constant,
we assume that $a^{(0)}_{\i0} = -1$. This condition on $\i0$ is considered to be valid hereinafter.
Let us build a parameterization of $\calD_r$ in the vicinity of $\tsS_0$; parameterization depends on the index $\i0$.
Note that we can not construct a global parameterization, since for any index $\i0$ there exists a point $\bfa$ of $\calD_r$ such that the $\i0$-th element of $\bfa$ is zero.

In the case of a series governed by an ordinary LRR($\bfa$), $\bfa\in \spR^{r+1}$, since the last coordinate of $\bfa$ is equal to $-1$,
the series is uniquely determined by the first $r$ elements of $\bfa$ and $r$  initial values of the series.
Then, applying the LRR to the initial data, which are taken from the series that is governed by the LRR,
we restore this series.
In the case of an arbitrary series in $\calD_r$, the approach is similar but a bit more complicated. For example,
we should take the boundary data ($\i0-1$ values at the beginning, and $r+1-\i0$ values at the end) instead of the $r$ initial values at the beginning of the series.

Denote $\calI(\i0) = \{1,\ldots, N\} \setminus \{\i0,\ldots, N-r-1+\i0\}$ and
$\calK(\i0) = \{1,\ldots,r+1\} \setminus \{\i0\}$ two sets of size $r$.
The set $\calI(\i0)$ consists of the indices of the series values (we call them boundary data), which
are enough to find all the series values with the help of the vector $\bfa_{\vecrow{\calK({\i0})} }\in \spR^{r}$
consisting of the elements of a vector $\bfa \in \spR^{r+1}$
with indices in $\calK({\i0})$; that is, $\bfa_{\vecrow{\calK({\i0})} }=(a_1,\ldots,a_{\i0-1},a_{\i0+1},a_{r+1})^\rmT$.

To simplify notation, let us introduce the operator $\fullop$: $\spR^{r} \to \spR^{r+1}$, which acts as follows. Let $\Ai0\in \spR^{r}$ and $\fullop(\Ai0) = \bfa$. Then $\bfa = (a_1, \ldots, $ $a_{r+1})^\rmT$ is such that $\bfa_{\vecrow{\calK({\i0})}} = \Ai0$ and $a_{\i0} = -1$; that is, $\Ai0 \in \spR^r$ is extended to $\bfa \in \spR^{r+1}$ by inserting $-1$ in the $\i0$-th position.
In this notation, $\bfa_{\vecrow{\calK({\i0})}} = \fullop^{-1}(\bfa)$.

Let $\bfZ=\bfZ(\bfa)$ be a matrix consisting of basis vectors of $\calZ(\bfa)$ and $\bfG(\Ai0) = \bfZ \left(\bfZ_{\row{\calI({\i0})}}\right)^{-1}$; here $\bfZ_{\row{\calI({\i0})}}\in \spR^{r\times r}$, $\bfG(\Ai0) \in \spR^{N\times r}$.
It follows from the explicit form of the parameterization given in \cite{Zvonarev.Golyandina2021} that
the parameterization mapping $S_{\tau}: \spR^{2r} \to \calD_r$ has the form
 		\begin{equation}\label{eq:param}
		\tsS = S_\tau(\Si0, \Ai0) = \bfG(\Ai0) \Si0.
		\end{equation}
Note that $\bfG(\Ai0)$ does not depend on the choice of basis vectors of $\calZ(\bfa)$ used for forming $\bfZ$. The expression $\bfG(\Ai0) \Si0$ gives the linear combination of the chosen basis vectors such that $\Si0 = (\tsS)_{\vecrow{\calI(\i0)}}$. Thus, $\tsS \in \calZ(\bfA)$ and has boundary data $\Si0$.

Note that for different series $\tsS_0\in \calD_r$ we may have different parameterizations of $\calD_r$ in vicinities of $\tsS_0$.
Moreover, for a fixed $\tsS_0$, there is a variety of parameterizations provided by different choices of the index $\i0$.

\section{Optimization}
\label{sec:optim}

Let us consider different numerical methods for solving the problem \eqref{eq:wls}.
First, note that we  search for a local minimum. Then, since the objective function is smooth in the considered
parameterization, one can apply the conventional weighted version of the Gauss-Newton method (GN), see \cite{nocedal2006numerical} for details.
However, this approach appears to be numerically unstable and has a high computational cost.

In \cite{Usevich2014}, the variable-projection method (VP) is used
for solving the minimization problem. When the reduced minimization problem
is solved again by the Gauss-Newton method, we will refer to it as VPGN.

We propose a similar (but different) approach called Modified Gauss-Newton method (MGN),
which appears to have some  advantages in comparison with VPGN that is one of the best methods
for solving the problem \eqref{eq:wls}. Below we will show that the MGN algorithm consists of different numerical sub-problems to be solved, which are more well-conditioned  than in the VPGN case (see \cite{Deuflhard.Hohmann2003}, where different properties of problems such as stability and well-conditioning are discussed); thereby, MGN allows a better numerically stable implementation.

Note that the considered methods are used for solving a weighted least-squares problem and therefore
we consider their weighted versions, omitting `weighted' in the names of the methods.

\medskip
Let us introduce the notation, which is used in this section.
For a matrix $\bfF = \spR^{N\times p}$, define its weighted pseudoinverse \cite{Stewart1989}
$\winverse{\bfF}{\bfW} = (\bfF^\rmT \bfW \bfF)^{-1}\bfF^\rmT \bfW$;
this pseudoinverse arises in the solution of the weighted linear least-squares problem $\min_\bfp \|\bfy - \bfF \bfp\|_{\bfW}^2$
with $\bfy \in \spR^{N}$, since its solution is equal to $\bfp_{\mathrm{min}} = \winverse{\bfF}{\bfW} \bfy$. In the particular case $\bfW = \bfI_N$, $\winverse{\bfF}{\bfW}$ is the ordinary pseudoinverse; we will denote it $\inverse{\bfF}$.
Denote the projection (it is oblique if $\bfW$ is not the identity matrix) onto the column space $\calF$ of a matrix $\bfF$ as $\Proj_{\bfF, \bfW} = \bfF \winverse{\bfF}{\bfW}$.
If it is not important which particular basis of $\calF$ is considered, we use the notation $\Proj_{\calF, \bfW}$.

\begin{remark}
\label{rem:complexF}
If the matrix $\bfF$ is complex, the above formulas and considerations are still valid with the change of the transpose $\bfF^\rmT$ to the complex conjugate $\bfF^\mathrm{H}$.
\end{remark}

The structure of this section is as follows.
After a brief discussion of the problem \eqref{eq:wls}, we start with the description of the methods GN and VP for a general optimization problem; then we apply these methods to \eqref{eq:wls} and finally present the new method MGN.

\subsection{Approaches for solving a general nonlinear least-squares problem}\label{sec:opt_gen}
	Let $\bfx \in \spR^N$ be a given vector and consider a general WLS minimization problem
	\begin{equation}
    \label{eq:gen_optim}
	\bfp^\star = \argmin_\bfp \|\bfx - S(\bfp)\|_{\bfW}^2,
	\end{equation}
	where $\bfp \in \spR^p$ is the vector of parameters, $S: \spR^p \to \spR^N$ is some parameterization of a subset of $\spR^N$ such that $S(\bfp)$ is a differentiable vector-function of $\bfp$,
$\bfW \in \spR^{N\times N}$ is
   a positive (semi-)definite symmetric matrix.

   If the problem \eqref{eq:gen_optim} is nonlinear, iterative methods with linearization
at each iteration are commonly used, such as the Gauss-Newton method or its variations \cite{nocedal2006numerical}.
One of the commonly used variations is the Levenberg-Marquardt method,
which is a regularized version of the Gauss-Newton method. This regularization improves the method
far from the minimum value and does not affect near the minimum. Therefore, in the paper, we consider the Gauss-Newton method without regularization.
We use the weighted Gauss-Newton method, which is a straightforward extension of the unweighted version with $\bfW=\bfI_N$.

\subsubsection{Gauss-Newton method} \label{sec:optim_wgn}

One iteration of the Gauss-Newton algorithm  with step $\gamma$ is
\begin{equation}
\label{eq:GN_P}
\bfp_{k+1} = \bfp_{k} + \gamma \winverse{\bfJ_{S}(\bfp_k)}{\bfW} (\bfx - S(\bfp_k)),
\end{equation}
where $\bfJ_{S}(\bfp_k)$ is the Jacobian matrix of $S(\bfp)$ at $\bfp_k$. Note that the iteration step \eqref{eq:GN_P} is uniquely defined for any positive semi-definite matrix, see Remark~\ref{rem:degenerate}.
The choice of step $\gamma$ is a separate problem. For example, one can apply the backtracking line search starting from $\gamma = 1$
and then decreasing the step if the next value is worse (that is, if the value of the objective functional increases).


An additional aim of the WLS problem is to find the approximation $S(\bfp^\star)$ of $\bfx$, where $\bfp^\star$ is the solution of \eqref{eq:gen_optim}.
 Then we can write \eqref{eq:GN_P} in the form of iterations of approximations:
\begin{equation} \label{eq:gnwls}
S(\bfp_{k+1}) = S \left(\bfp_{k} + \gamma \winverse{\bfJ_{S}(\bfp_k)}{\bfW} \left(\bfx - S(\bfp_k)\right)\right).
\end{equation}

The following remark explains the approach, which underlies the Modified Gauss-Newton method proposed in this paper (see Section~\ref{sec:MGN_opt}).
\begin{remark} \label{rem:replacement}
	The iteration step \eqref{eq:gnwls} can be varied by means of the change of $S(\bfp_{k+1})$ to $\widetilde S(\bfp_{k+1})$, where $\widetilde S(\bfp_{k+1})$ is such that  $\|\bfx - \widetilde S(\bfp_{k+1})\|_\bfW \le \|\bfx - S(\bfp_{k+1})\|_\bfW$. This trick is reasonable if $\widetilde S(\bfp_{k+1})$ can be calculated
faster and/or in a more stable way than $S(\bfp_{k+1})$. 
\end{remark}

\subsubsection{Variable projection}
	Let $\bfp = \left(\begin{matrix}\bfb\\\bfc\end{matrix}\right)\in \spR^p$, $\bfb\in \spR^{p_1}$, $\bfc\in \spR^{p_2}$.
Consider the (weighted) least-squares problem \eqref{eq:gen_optim}, where $S(\bfp)$ is linear in $\bfc$ and the nonlinear part is defined by
$\bfG(\bfb)\in \spR^{N\times p_2}$:
\bea
S(\bfp) = \bfG(\bfb) \bfc.
\eea
This problem can be considered as a problem of projecting the data vector $\bfx$ onto a given set:
	\begin{equation}
\label{eq:full_optim}
	\min_{\bfy \in \calD} \|\bfx - \bfy \|_{\bfW}, \quad \text{where} \quad \calD =
    \Big\{\bfG(\bfb) \bfc \mid \left(\begin{matrix}\bfb\\\bfc\end{matrix}\right)\in \spR^p\Big\}.
	\end{equation}
	 Here $\{\varphi(z)\mid z\in \calC\}$ means the set of values of $\varphi(z)$ for $z\in \calC$.
    The variable projection method takes advantage of the known explicit solution of the subproblem:
	\begin{equation*}
	C^\star(\bfb) = \argmin_{\bfc} \|\bfx - \bfG(\bfb) \bfc \|_{\bfW} = \winverse{\bfG(\bfb)}{\bfW} \bfx.
	\end{equation*}

\smallskip	
	Denote $S^\star(\bfb) = \bfG(\bfb) C^\star(\bfb)$, $\calG(\bfb) = \{\bfG(\bfb)\bfc \mid \bfc\in \spR^{p_2} \}$. Then
	\begin{equation} \label{eq:minZ}
	S^\star(\bfb) = \argmin_{\bfs \in \calG(\bfb)} \|\bfx - \bfs \|_{\bfW} = \Proj_{\calG(\bfb), \bfW}\bfx.
	\end{equation}
	Thus, we can reduce the problem \eqref{eq:full_optim} to the projection onto a subset $\calD^\star \subset \calD$ and
thereby to the optimization in the nonlinear part of parameters only:
	\begin{equation}\label{eq:vpprinciple}
	\min_{\bfy \in \calD^\star} \|\bfx - \bfy \|_{\bfW} \quad \text{with} \quad \calD^\star = \{S^\star(\bfb) \mid \bfb\in\spR^{p_1}\}.
	\end{equation}
This is called ``variable projection'' principle (see \cite{Golub.Pereyr2003} for the case of the Euclidean norm).

\subsection{Known iterative methods for low-rank approximation}
Let us turn from a general nonlinear least-squares problem \eqref{eq:gen_optim} to the specific problem \eqref{eq:wls}.

A variation from the standard way of the use of iterative methods is that the parameterization $S_{\tau}(\bfp)$, $\bfp=(\Si0, \Ai0)$ (which is based on $\i0$)
is changed at each iteration in a particular way.
At $(k+1)$-th iteration, the parameterization is constructed in the vicinity of  $\bfa_0 = \bfa^{(k)}$. The index $\i0$, which determines the parameterization, is chosen in such a way to satisfy $a^{(0)}_\i0 \neq 0$.
We propose the following approach to the choice of $\i0$.
Let $\i0$ be the index of the maximum absolute entry of $\bfa_0$. Since the parameterization is invariant to the multiplication of $\bfa_0$ by a constant, it can be assumed that $a^{(0)}_\i0 = -1$ and $|a^{(0)}_i| \le 1$ for any $i$, $1 \le i \le r+1$.

\subsubsection{Weighted Gauss-Newton method}

The Gauss-Newton algorithm can be applied to the problem \eqref{eq:wls} in a straightforward manner, taking into consideration that the parameterization $S_{\tau}$ may be changed at each iteration.
The Gauss-Newton iteration has the form $\bfp_{k+1} = \bfp_k + \gamma \winverse{\bfJ_{S_{\tau}}(\bfp_k)}{\bfW} (\tsX - S_{\tau}(\bfp_k))$.

To apply the method, $S_{\tau}(\bfp_k)$ and the Jacobian matrix $\bfJ_{S_{\tau}}(\bfp_k)$
should be calculated. Formally, their computing can be implemented; however, the direct calculation is not
numerically stable and very time-consuming.

\subsubsection{Gauss-Newton method with variable projection}
\label{sec:VPGN}
The explicit form of the parameterization $S_{\tau}(\bfp) = S_{\tau}(\Si0, \Ai0)$ given in \eqref{eq:param}, where $\Si0$ is presented in $S_{\tau}(\Si0, \Ai0)$ in a linear manner,
allows one to apply the variable projection principle.

Assume that $\tsS_0$ is governed by a GLRR($\bfa_0$) with $a^{(0)}_{\i0} = -1$ and consider the problem \eqref{eq:wls} in the vicinity of the series $\tsS_0 \in \calD_r$.

Substitute in \eqref{eq:vpprinciple} $\calD = \overline \calD_r$, $\calD^\star = \calD_r^\star \subset \overline \calD_r$, where $\calD_r^\star = \{\Proj_{\calZ(\fullop(\Ai0)), \bfW}\tsX \mid \Ai0 \in \spR^{r}\}$, $\bfb = \Ai0$, $\bfG(\bfb) = \bfG$, where $\bfG = \bfZ \left(\bfZ_{\row{\calI({\i0})}}\right)^{-1}$ (see \eqref{eq:param}), $C^\star(\bfb) = \winverse{\bfG}{\bfW} \tsX$. Then
\begin{equation}
S_{\tau}^\star(\Ai0) = \Proj_{\calZ(\bfa), \bfW}(\tsX)
\end{equation}
where $\bfa = \fullop(\Ai0)$, and
 we obtain the equivalent problem for projecting the elements of the set $\overline \calD_r$ to the subset $\calD_r^\star$,
where the parameter $\Si0$ is eliminated:
\begin{equation}\label{eq:wlsvp_set}
\tsY^\star = \argmin_{\tsY \in \calD_r^\star} \| \tsX - \tsY \|_{\bfW}.
\end{equation}
or, the same,
\begin{equation}\label{eq:wlsvp1}
\Ai0^\star = \argmin_{\Ai0 \in \spR^r} \| \tsX - S_{\tau}^\star(\Ai0) \|_{\bfW},
\end{equation}
Thus, for the numerical solution of the equation \eqref{eq:wls}, it is sufficient to consider iterations for the nonlinear part of the parameters.
This is the VP approach used in \cite{Usevich2012, Usevich2014}.

Let us denote $\bfJ_{S_{\tau}^\star}(\Ai0)$ the Jacobian matrix of $S_{\tau}^\star(\Ai0)$.
Then the iterations of the Gauss-Newton method for solving the problem \eqref{eq:wlsvp1} have the form
\begin{equation}
\label{eq:gauss_simple}
\Ai0^{(k+1)} = \Ai0^{(k)} + \gamma \winverse{\bfJ_{S_{\tau}^\star}(\Ai0^{(k)})}{\bfW} \big(\tsX - S_{\tau}^\star(\Ai0^{(k)})\big).
\end{equation}
The VPGN algorithm together with an explicit form of $\bfJ_{S_{\tau}^\star}(\Ai0^{(k)})$ is presented in \cite{Zvonarev.Golyandina2021}; its iteration step is described in Algorithm~\ref{alg:gauss_newton_vp_step}.

\subsection{Modified Gauss-Newton method for low-rank approximation}
\label{sec:MGN_opt}
In this section, we propose a new iterative method for the problem \eqref{eq:wls}, which is a modified Gauss-Newton method.

Let us return to the problem with the full set of parameters $(\Si0,\Ai0)$ and apply the approach that is described in Remark~\ref{rem:replacement},
 with $\widetilde S(\bfp) = S_{\tau}^\star(\Ai0)$. We can do it, since $S_{\tau}^\star(\Ai0) = \Proj_{\calZ(\fullop(\Ai0)), \bfW}\,\tsX$ and therefore \eqref{eq:minZ} is valid with $\calG(\Ai0) = \calZ(\fullop(\Ai0))$.
 Thus, we can consider $S_{\tau}^\star\big( \Ai0^{(k+1)}\big) \in \calD_r^\star$ as the result of the $(k+1)$-th iteration instead of
$S_{\tau}\big(\Si0^{(k+1)}, \Ai0^{(k+1)}\big) \in \overline \calD_r$.
It appears (see \cite{Zvonarev.Golyandina2021}) that then we can use more stable numerical calculations for the iteration implementation.
The proposed modification is similar to variable projections, since we can omit the part $\Si0$ of parameters.

Thus, we introduce the MGN iteration in the form
\begin{equation}
\label{eq:iterGNA}
\Ai0^{(k+1)}
= \Ai0^{(k)} + \gamma \left( \winverse{\bfJ_{S_{\tau}}(\Si0^{(k)}, \Ai0^{(k)})}{\bfW} \big(\tsX - S_{\tau}^\star(\Ai0^{(k)})\big) \right) _{\col{\{r+1, \ldots, 2r\}}},
\end{equation}
where $\Si0^{(k)}$ are the corresponding boundary data taken from $S_\tau^\star(\Ai0^{(k)})$, i.e. $\Si0^{(k)} = \left( S_\tau^\star(\Ai0^{(k)}) \right)_{\calI(\tau)}$.
As well as in the variable projection method with the iteration step \eqref{eq:gauss_simple}, $S_\tau^\star(\Ai0^{(k+1)}) \in \calD_r^\star$ for each $k$.

\begin{theorem}\label{th:equivalency}
Let $\bfW^{1/2} \bfJ_{S_{\tau}}(\Si0^{(k)}, \Ai0^{(k)})$ have full rank.
Denote $\bfS = \calT_{r+1} \left(\Proj_{\calZ(\fullop(\bfa^{(k)})), \bfW}\,\tsX \right)$,
$\bfM = - \left(\bfS_{\row{\calK({\i0})}}\right)^\rmT$.
Then the iteration step \eqref{eq:iterGNA} is equivalent to
\begin{equation}
\label{eq:iterGNfinal}
	\Ai0^{(k+1)}
	= \Ai0^{(k)} + \gamma \winverse{\big(\bfI_N - \Proj_{\calZ(\fullop(\Ai0^{(k)})), \bfW}\big)\widehat \bfF_{\bfa}}{\bfW} \big(\tsX - S_{\tau}^\star(\Ai0^{(k)})\big),
\end{equation}
where $\widehat \bfF_{\bfa} \in \spC^{N \times 2r}$ is an arbitrary matrix satisfying the equality $\bfQ^\rmT(\fullop(\Ai0^{(k)})) \widehat \bfF_{\bfa} = \bfM$.
\end{theorem}
Let us fix the iteration number $k$.
Denote by $\bfF_\bfs = \left(\bfJ_{S_\tau} \right)_{\col{\{1, \ldots, r\}}}$ the first $r$ columns of the Jacobian matrix $\bfJ_{S_\tau} = \bfJ_{S}(\Si0^{(k)}, \Ai0^{(k)})$, and by $\bfF_\bfa = \left(\bfJ_{S_\tau} \right)_{\col{\{r+1, \ldots, 2r\}}}$ the last $r$ columns of $\bfJ_{S_\tau}$.

Before proving the theorem, let us write down the statements of several propositions of \cite{Zvonarev.Golyandina2021}.

By definition, the tangent space at the point $\tsS$ coincides with $\colspace\left(\bfJ_{S_{\tau}}(\Si0, \Ai0)\right)$.
Note that the tangent space is invariant with respect to the choice of a certain parameterization of $\calD_r$ in the vicinity of $\tsS$.

Define by $\bfa^2$ the acyclic convolution of $\bfa$ with itself:
\begin{equation*}
	\bfa^2 = (a^{(2)}_i) \in \spR^{2r+1}, \quad a^{(2)}_i = \sum_{j=\max(1, i - r)}^{\min(i, r+1)} a_j a_{i - j + 1}.
\end{equation*}
		
    \begin{lemma}[{\cite{Zvonarev.Golyandina2021}}]
    	\label{eqa:derivS}
    	$\bfQ^\rmT(\bfa) \bfF_\bfs=\bfzero_{(N-r)\times r}$; $\colspace(\bfF_\bfs) = \calZ(\bfa)$.
    \end{lemma}	

    \begin{lemma}[{\cite{Zvonarev.Golyandina2021}}]
    	\label{eqa:derivA}
    	$\bfQ^\rmT(\bfa) \bfF_{\bfa} =  - (\bfS_{\row{\calK({\i0})}})^\rmT$, where $\bfS = \calT_{r+1}(\tsS)$;
    	$\colspace(\bfF_\bfa) \subset \calZ(\bfa^2)$.
    \end{lemma}

\begin{theorem}[{\cite{Zvonarev.Golyandina2021}}]
\label{th:tangent}
	The tangent space of $\calD_r$ at the point $\tsS$ has dimension $2r$ and is equal to $\calZ(\bfa^2)$.
\end{theorem}

\begin{proof}[Proof of Theorem~\ref{th:equivalency}]
Let us rewrite the weighted pseudoinverse in the \eqref{eq:iterGNA} as
\begin{equation*}
	\left( \winverse{\bfJ_{S_{\tau}}(\Si0^{(k)}, \Ai0^{(k)})}{\bfW} \big(\tsX - S_{\tau}^\star(\Ai0^{(k)})\big) \right) _{\col{\{r+1, \ldots, 2r\}}}
	= \left( \inverse{\left(\bfW^{1/2} \bfJ_{S_{\tau}}(\Si0^{(k)}, \Ai0^{(k)})\right)} \bfW^{1/2}\big(\tsX - S_{\tau}^\star(\Ai0^{(k)})\big) \right) _{\col{\{r+1, \ldots, 2r\}}}.
\end{equation*}	
  Applying the Frisch-Waugh-Lovell theorem \cite{Lovell2008simple} about the partitioned regression to the obtained pseudoinverse for regressors $\bfW^{1/2} \bfF_{\bfs}$ and $\bfW^{1/2} \bfF_{\bfa}$, we get the following sequence of equalities:
\begin{multline*}
\left( \inverse{\left(\bfW^{1/2} \bfJ_{S_{\tau}}(\Si0^{(k)}, \Ai0^{(k)})\right)} \bfW^{1/2}\big(\tsX - S_{\tau}^\star(\Ai0^{(k)})\big) \right) _{\col{\{r+1, \ldots, 2r\}}} = \\
\inverse{\left( (\bfI_N - \Proj_{\bfW^{1/2} \bfF_{\bfs}})\bfW^{1/2} \bfF_{\bfa} \right)}(\bfI_N - \Proj_{\bfW^{1/2} \bfF_{\bfs}}) \bfW^{1/2} (\tsX - S_{\tau}^\star(\Ai0^{(k)})) = \\
\winverse{(\bfI_N - \Proj_{\calZ(\fullop(\Ai0^{(k)})), \bfW})\bfF_{\bfa}}{\bfW}(\bfI_N - \Proj_{\calZ(\fullop(\Ai0^{(k)})), \bfW}) (\tsX - \Proj_{\calZ(\fullop(\Ai0^{(k)})), \bfW}\,\tsX).
\end{multline*}	
Since $\bfI_N - \Proj_{\calZ(\fullop(\Ai0^{(k)})), \bfW}$ is a projector, $(\bfI_N - \Proj_{\calZ(\fullop(\Ai0^{(k)})), \bfW})^2 = \bfI_N - \Proj_{\fullop(\calZ(\bfa^{(k)})), \bfW}$. Thus, we obtain the following iteration equivalent to \eqref{eq:iterGNA}:
\begin{equation}
\Ai0^{(k+1)}
= \Ai0^{(k)} + \gamma \winverse{(\bfI_N - \Proj_{\calZ(\fullop(\Ai0^{(k)})), \bfW})\bfF_{\bfa}}{\bfW}(\bfI_N - \Proj_{\calZ(\fullop(\Ai0^{(k)})), \bfW})\,\tsX.
\end{equation}
By Lemma~\ref{eqa:derivA}, $\bfQ^\rmT(\fullop(\Ai0^{(k)})) \bfF_\bfa = \bfM$. By the theorem's conditions, $\bfQ^\rmT(\fullop(\Ai0^{(k)})) \widehat \bfF_\bfa = \bfM$. Thus, $\bfQ^\rmT(\fullop(\Ai0^{(k)}))(\bfF_\bfa - \widehat \bfF_\bfa) = \bfzero_{(N-r)\times r}$. Since $\calQ(\bfa)$ is the orthogonal complement to $\calZ(\bfa)$, $(\bfI_N - \Proj_{\calZ(\fullop(\Ai0^{(k)})), \bfW})(\bfF_\bfa - \widehat \bfF_\bfa) = \bfzero_{N\times r}$, which finishes the proof.
\end{proof}

\begin{remark}
In the case when $S_{\tau}(\Si0^{(k)}, \Ai0^{(k)}) \in \calD_r$, the Jacobian matrix $\bfJ_{S_{\tau}}(\Si0^{(k)}, \Ai0^{(k)})$ has full rank, according to Theorem~\ref{th:tangent}. This is sufficient for validity of the condition of Theorem~\ref{th:equivalency} if $\bfW$ has full rank. In the case of a rank-deficient matrix $\bfW$, the condition of Theorem~\ref{th:equivalency} is discussed in Remark~\ref{rem:degenerate} with $\bfF = \bfJ_{S_{\tau}}(\Si0^{(k)}, \Ai0^{(k)})$.
\end{remark}

Thus, we have constructed the equivalent version \eqref{eq:iterGNfinal} of the iteration step \eqref{eq:iterGNA} in such a way to reduce its complexity to the computational costs of computing the projections to $\calZ(\fullop(\Ai0)) = \calZ(\bfa)$ and calculating the matrices $\widehat \bfF_{\bfa}$ for different $\bfa$.
A numerically robust algorithm for calculating the iteration step \eqref{eq:iterGNfinal} is given in Algorithm~\ref{alg:gauss_newton_our_step}.

\section{Algorithms of the VPGN and MGN methods}
\label{sec:MGNand VPGN}

The common scheme of the VPGN and MGN algorithms is given in Algorithm~\ref{alg:common_scheme}.
\begin{algorithm}
	\caption{The common scheme of VPGN and MGN}
	\label{alg:common_scheme}
    \Input{$\tsX \in \spR^N$, $\bfa_0\in \spR^{r+1}$, a stopping criterion STOP.}
	\begin{algorithmic}[1]
		\State{Set $k = 0$,  $\bfb^{(0)} = \bfa_0$.}
		\Repeat{}
		\State{Choose $\i0$ such that $b_\i0^{(k)}\neq 0$; for example, find $\i0 = \argmax_{i} |b_i^{(k)}|$. Calculate $\bfa^{(k)} = -\bfb^{(k)}/b_\tau^{(k)}$ to obtain $a_\tau^{(k)} = -1$
and take $\Ai0^{(k)} = \fullop^{-1}(\bfa^{(k)})$.}
		\State{Calculate $\tsS_k  = S_\tau^\star(\Ai0^{(k)})$, where $S_\tau^\star(\Ai0^{(k)})  = \Proj_{\calZ(\bfa^{(k)}) , \bfW}\tsX$.}
		\State{Calculate the direction $\Delta_k = \Delta(\Ai0^{(k)}, \tsS_k)$.}
		\State{Choose the size $\gamma_k$ and perform the step of size $\gamma_k$ in the descent direction given by $\Delta_k$. For example, find $\gamma_k =\gamma$, $0\le \gamma\le 1$, such that
          $\|\tsX - S_{\tau}^\star(\Ai0^{(k)} + \gamma \Delta_k) \|_{\bfW} \le \|\tsX - S_\tau^\star(\Ai0^{(k)}) \|_{\bfW}$
          by the backtracking line search method \cite[Section 3.1]{nocedal2006numerical}.}
		\State{Set $\Ai0^{(k+1)} = \Ai0^{(k)} + \gamma_k \Delta_k$, $\bfb^{(k+1)} = \fullop(\Ai0^{(k+1)})$.}
		\State{Set $k = k+1$.}
		\Until{STOP}
		\State\Return{$\widetilde \tsS = S_{\i0}^\star(\Ai0^{(k)})$ as an estimate of the signal.}
	\end{algorithmic}
\end{algorithm}

As it is shown in \cite{Zvonarev.Golyandina2021}, the calculation of the direction in step 5 is the ``bottleneck'' of the VPGN algorithm, since the calculation of $S_{\tau}^\star$ in steps 4 and 6 can be stably implemented (see S-VPGN in \cite{Zvonarev.Golyandina2021}). Moreover, the computational cost of step 5 is very large for the case of banded weight matrices $\bfW$, which correspond to autoregressive noise. The MGN method is designed to improve step 5.

Thus, we have Algorithm~\ref{alg:common_scheme} with

\begin{eqnarray}
\label{eq:delta_VPGN}
\text{VGPN:\quad } \Delta(\Ai0^{(k)}, \tsS_k) &=& \winverse{\bfJ_{S_{\tau}^\star}(\Ai0^{(k)})}{\bfW} \big(\tsX - \tsS_k\big),\\
\label{eq:delta_MGN}
\text{MGN:\quad } \Delta(\Ai0^{(k)}, \tsS_k) &=& \winverse{\big(\bfI_N - \Proj_{\calZ(\fullop(\Ai0^{(k)})), \bfW}\big)\widehat \bfF_{\bfa}}{\bfW} \big(\tsX - \tsS_k\big),
\end{eqnarray}
where $\tsS_k = S_{\tau}^\star(\Ai0^{(k)}) = \Proj_{\calZ(\fullop(\Ai0^{(k)})), \bfW}\,\tsX$.

The stable versions of calculating the projection are described in \cite{Zvonarev.Golyandina2021} (the algorithms S-VPGN and S-VPGN-H).
All the version of VPGN still have a ``bottleneck'' in calculation of $\bfJ_{S_{\tau}^\star}(\Ai0^{(k)})$, while the whole MGN algorithm can be stably implemented.

\subsection{Stable auxiliary algorithms}

Algorithm~\ref{alg:fourier_basis_A} below contains the stable implementation of computing the basis for calculating the projection $\Proj_{\calZ(\fullop(\Ai0^{(k)})),\bfW}$ \cite{Zvonarev.Golyandina2021}.
The step of MGN additionally contains the calculation of $\widehat \bfF_{\bfa}$. Applying the technique of \cite{Zvonarev.Golyandina2021}, we can derive Algorithm~\ref{alg:fourier_grad} for calculating $\widehat \bfF_{\bfa}$.

Let us briefly describe the approach of \cite{Zvonarev.Golyandina2021}.
Denote
 \be
 \label{eq:pol_z}
 g_{\bfa}(z) = \sum_{k=0}^{r} a_{k+1}z^k
 \ee
   the complex polynomial with coefficients $\bfa = (a_1,\ldots,a_{r+1})^\mathrm{T}$; we do not assume that the leading coefficient is non-zero.

Let the circulant matrix $\bfC(\bfa)$ be the extension of the partial circulant $\bfQ^\rmT(\bfa)$.
Then the construction of a basis of $\calZ(\bfa)$ can be reduced to
solving the system of linear equations
\be
\label{eq:lineqMGN}
    \bfC (\bfa) \bfv_k = \bfe_{N-k+1},\ k=1,\ldots,r.
\ee
The eigenvalues of $\bfC(\bfa)$ coincide with the values of the polynomial $g_\bfa(z)$ in nodes of the equidistant grid $\calW = \left\{\exp\big(\frac{\unit 2 \pi j}{N}\big), \; j = 0, \ldots, N-1\right\}$ on the complex unit circle $\spT= \{z \in \spC : |z| = 1 \}$. Therefore, the nondegeneracy of $\bfC(\bfa)$ is equivalent to that there are no roots of the polynomial $g_\bfa(z)$ in $\calW$.

Let us define the unitary matrix
\begin{equation}
\label{eq:eqi_grid}
\bfT_M(\alpha) = \diag\left((1, e^{\unit \alpha}, \ldots, e^{\unit (M-1) \alpha})^\rmT\right),
\end{equation}
where $\alpha$ is a real number, $M$ is a natural number.

\begin{remark}
\label{rem:shifting}
In the exact arithmetic, an arbitrary small non-zero value of the smallest eigenvalue of a matrix provides its non-degeneracy. However, in practice, the numerical stability and accuracy of matrix calculations depend on the condition numbers of matrices.
Therefore, the aim of the choice of a proper $\alpha$ is to do the condition number of $\bfC(\tilde \bfa(\alpha))$, where $\tilde \bfa = \tilde \bfa(\alpha) = \left(\bfT_{r+1}(-\alpha)\right) \bfa$, as small as possible.
This minimization problem can be approximately reduced to the problem
of maximization of the smallest eigenvalue $|\lambda_\text{min}(\alpha)| = \min_{z \in \calW(\alpha)} | g_\bfa(z) |$ of $\bfC(\tilde \bfa(\alpha))$, since the maximal eigenvalue is not larger than $\max_{z \in \spT} |g_\bfa(z)|$.
\end{remark}

Denote $\calF_N$ and $\calF^{-1}_N$ the Fourier transform and the inverse Fourier transform for series of length $N$, respectively.
Define $\calF_N(\bfX) = [\calF_N(\bfx_1): \ldots: \calF_N(\bfx_r)]$, where $\bfX=[\bfx_1 : \ldots : \bfx_r]$; the same for $\calF_N^{-1}(\bfY)$. Algorithm~\ref{alg:fourier_basis_A} exploits the fact that the operations with circulant matrices can be implemented by means of fast Fourier transform \cite{Korobeynikov2010}.

\begin{algorithm}
	\caption{Calculation of a basis of $\calZ(\bfa)\subset \spC^N$ \cite{Zvonarev.Golyandina2021}}
	\label{alg:fourier_basis_A}
    \Input{$\bfa \in \spR^r$.}
	\begin{algorithmic}[1]
		\State Find $\alpha_0 = \argmax_{-\pi/N \le \alpha < \pi/N} \min_{z \in \calW(\alpha)} | g_\bfa(z) |$ by means of a 1D numerical optimization method.
		\State Calculate  the vector $\bfa_g = (a_{g, 0}, \ldots, a_{g, N-1})^\rmT$ consisting of the eigenvalues of $\bfC(\widetilde \bfa)$ by $a_{g, j} = g_\bfa\big(\exp(\unit (\frac{2 \pi j}{N} - \alpha_0)\big)$, $j = 0, \ldots, N-1$, where $\widetilde \bfa = \left(\bfT_{r+1}(-\alpha)\right) \bfa$; $\bfA_g = \diag(\bfa_g)$.
		\State Calculate the matrices $\bfR_r = \calF_N([\bfe_{N-r+1}: \ldots: \bfe_N])$ and  $\bfL_r = \bfA_g^{-1} \bfR_r$.
		\State Compute $\bfU_r$: find $\bfO_r$ such that $\bfL_r \bfO_r$ consists of orthonormal columns ($\bfO_r$ can be found by either the QR factorization or the SVD);
calculate $\bfB=\bfR_r \bfO_r$ and $\bfU_r = \bfA_g^{-1} \bfB$.
		\State Compute $\widetilde \bfZ = \calF_N^{-1}(\bfU_r)$.
		\State \Return $\bfZ = (\bfT_{N}(-\alpha_0)) \widetilde \bfZ \in \spC^{N\times r}$, whose columns form an orthonormal basis of $\calZ(\bfa)$, $\alpha_0$ and $\bfA_g$.
	\end{algorithmic}
\end{algorithm}

 Let us now turn to calculating the matrix $\widehat \bfF_{\bfa}$. According to Theorem~\ref{th:equivalency}, it is sufficient to find an arbitrary matrix such that $\bfQ^\rmT(\bfa) \widehat \bfF_{\bfa} = \bfM$, where $\bfM \in \spR^{(N-r) \times r}$ is defined in Theorem~\ref{th:equivalency}. Therefore, it is sufficient to solve the following systems
of linear equations:
\be
\label{eq:lineqMGN2}
    \bfC (\bfa) \widehat \bfF_{\bfa} = \left(\begin{matrix}\bfM\\\bfzero_{r \times r}\end{matrix}\right).
\ee
The following lemma is a direct application of the theorem about the solution of a linear system of equations given by a circulant matrix \cite{Davis2012}.

\begin{lemma}
	\label{lemma:ev_circulant_pre}
 Define $\widehat \bfF_\bfa = \calF_N^{-1}(\bfA_g^{-1} \widehat \bfR_r)$, where $\widehat \bfR_r = \calF_N \left( \left(\begin{matrix}\bfM\\\bfzero_{r \times r}\end{matrix}\right) \right)$. Then $\bfQ^\rmT(\bfa) \widehat \bfF_\bfa = \bfM$, i.e. $\widehat \bfF_\bfa$ satisfies the conditions of Theorem \ref{th:equivalency}.
\end{lemma}

Remark~\ref{rem:shifting} together with Lemma~\ref{lemma:ev_circulant_pre} provide Algorithm~\ref{alg:fourier_grad}.

\begin{algorithm}
	\caption{Calculation of a matrix $\widehat \bfF_\bfa$ in \eqref{eq:iterGNfinal}}
	\label{alg:fourier_grad}
    \Input{$\bfa \in \spR^r$ and a series $\tsS\in \spR^N$ governed by the GLRR($\bfa$).}
	\begin{algorithmic}[1]
		\State{Compute $\alpha_0$, $\bfA_g$ using steps 1 and 2 of Algorithm \ref{alg:fourier_basis_A}. }
        \State{Construct $\bfM = - (\bfS_{\row{\calK({\i0})}})^\rmT$, where $\bfS = \calT_{r+1} \left( \tsS \right)$.}
		\State{Calculate $\widetilde \bfM = \left(\begin{matrix} (\bfT_{N-r}(\alpha_0)) \bfM\\\bfzero_{r \times r}\end{matrix}\right)$.}
		\State{Calculate $\widehat \bfR_{r} = \calF_N(\widetilde \bfM)$ and  $\widetilde \bfF_\bfa = \calF_N^{-1}(\bfA_g^{-1} \widehat \bfR_{r})$.}
		\State\Return{$\widehat \bfF_\bfa = (\bfT_{N}(-\alpha_0)) \widetilde \bfF_\bfa  \in \spC^{N\times 2r}$.}
	\end{algorithmic}
\end{algorithm}

\begin{remark}
\label{rem:hornerscheme}
In \cite{Zvonarev.Golyandina2021}, an enhanced version of Algorithm~\ref{alg:fourier_basis_A} using the compensated Horner scheme of \cite{Graillat2008} is considered. In the version with the compensated Horner scheme, this scheme is applied to calculating the values of the polynomial $g$ and multiplication $\bfR_r \bfO_r$. It is shown that this way of calculations considerably improves the stability of the projection algorithm.
The construction of Algorithm~\ref{alg:fourier_grad} allows one to improve its stability by the change of the first step
``Compute $\alpha_0$, $\bfA_g$ using Algorithm \ref{alg:fourier_basis_A}'' to its enhanced version.
\end{remark}

\subsection{Algorithms for the iteration step in VPGN and MGN}
Thus, both VPGN and MGN are particular cases of Algorithm~\ref{alg:common_scheme} with different directions of the line search given by formulas \ref{eq:delta_VPGN} and \ref{eq:delta_MGN} and presented in Algorithms~\ref{alg:gauss_newton_vp_step} and \ref{alg:gauss_newton_our_step}.

\begin{algorithm}
	\caption{Calculating the direction in VPGN}
	\label{alg:gauss_newton_vp_step}
    \Input{$\tsX \in \spR^N$, $\Ai0^{(k)}\in \spR^{r+1}$.}
	\begin{algorithmic}[1]
		\State{Calculate $\tsS_k  = S_\tau^\star(\Ai0^{(k)})$ as $\tsS_k  = \Proj_{\calZ(\bfa^{(k)}) , \bfW}\tsX$.}
		\State{Calculate $\bfJ_{S_\tau^\star}(\Ai0^{(k)})$.}
		\State{Calculate
			$\Delta_k = \winverse{\bfJ_{S_\tau^\star}(\Ai0^{(k)})}{\bfW} (\tsX - \tsS_k)$}.
		\State\Return{$\Delta_k$}
	\end{algorithmic}
\end{algorithm}

\begin{algorithm}
	\caption{Calculating the direction in MGN}
	\label{alg:gauss_newton_our_step}
    \Input{$\tsX \in \spR^N$, $\Ai0^{(k)}\in \spR^{r+1}$.}
	\begin{algorithmic}[1]
\State{Calculate $\tsS_k  = S_\tau^\star(\Ai0^{(k)})$  as $\tsS_k  = \Proj_{\calZ(\bfa^{(k)}) , \bfW}\tsX$.}
\State{Calculate $\widehat \bfF_{\bfa^{(k)}}$ by Algorithm \ref{alg:fourier_grad} with $\bfa = \bfa^{(k)}$
and $\tsS=\tsS_k$.}
\State{Calculate
$\Delta_k = \winverse{\left(\bfI_N - \Proj_{\calZ(\bfa^{(k)}), \bfW}\right)\widehat \bfF_{\bfa^{(k)}}}{\bfW} (\tsX - \tsS_k)$.}
		\State\Return{$\Delta_k$}
	\end{algorithmic}
\end{algorithm}

The technique for calculating the projector $\Proj_{\calZ(\bfa), \bfW}$ developed in \cite{Zvonarev.Golyandina2021} uses the basis consisting of the vectors found by Algorithm~\ref{alg:fourier_basis_A}. This technique supplemented by Algorithm~\ref{alg:fourier_grad} allows one to create fast and stable implementation of MGN (Algorithm~\ref{alg:gauss_newton_our_step}).
The calculation of $\bfJ_{S_\tau^\star}(\Ai0^{(k)})$ does not allow for both a fast and stable implementation of VPGN, even if Algorithm~\ref{alg:fourier_basis_A} is used for calculating the projector $\Proj_{\calZ(\bfa), \bfW}$ (the S-VPGN algorithm in \cite{Zvonarev.Golyandina2021}).

\subsection{The case of degenerate $(2p+1)$-diagonal $\bfW$}
\label{sec:missing}
Let us show how Algorithm~\ref{alg:gauss_newton_our_step} in MGN can be implemented for a time series with missing values. Let a symmetric positive-definite matrix $\bfW_0$ be given for the whole time series including missing entries. Then the weight matrix $\bfW$ is constructed from $\bfW_0$ by setting the values of its columns and rows with indices equal to the entries of missing values to zero.
In the case of AR($p$) noise, the matrix $\bfW_0$ (and therefore $\bfW$) is $(2p+1)$-diagonal.

Thus, consider a $(2p+1)$-diagonal matrix $\bfW$ which can be degenerate.

In the considered algorithms, there is a problem of calculating the weighted pseudoinverses and projections if the weight matrix $\bfW$ is degenerate.
Let us extend the definition of $\winverse{\bfF}{\bfW}$ and $\Proj_{\bfF, \bfW}$ for some matrix $\bfF$.

Consider a degenerate case when $\bfF^\rmT \bfW \bfF$ is not positive definite or, the same,
$\bfW ^{1/2}\bfF$ is rank-deficient ($\bfW ^{1/2}$ is the principal square root of $\bfW$). Then we can use a different representation for the weighted pseudoinverse:
$\winverse{\bfF}{\bfW} = \inverse{(\bfW^{1/2} \bfF)} \bfW^{1/2}$.
This corresponds to the minimum-(semi)norm solution of the corresponding WLS problem $\min_\bfp \|\bfy - \bfF \bfp\|_{\bfW}^2$.
Although
the projection $\Proj_{\bfF, \bfW}$ is generally not uniquely defined in the degenerate case, we will consider its uniquely defined version given by the formula  $\Proj_{\bfF, \bfW} = \bfF \winverse{\bfF}{\bfW}$.

\begin{remark}
\label{rem:degenerate}
Note that the matrix $\bfW ^{1/2}\bfF$ is rank-deficient if $\bfF$ is rank-deficient.
However, for a full-rank $\bfF$ and a degenerate $\bfW$, $\bfW ^{1/2}\bfF$ is not necessarily rank-deficient.
For example, if the orthogonal projections of the columns of $\bfF$ on $\colspace(\bfW)$ are linearly independent,
then  $\bfW ^{1/2}\bfF$ is full-rank.
\end{remark}

\paragraph{Calculation of weighted projection onto subspace with a given basis}
\label{sec:Pi_weighted}
     Consider the Cholesky decomposition $\Sigminus = \bfC^\rmT \bfC$; here $\bfC$ is an upper triangular matrix with $p$ nonzero superdiagonals \cite[p. 180]{GoVa13}.

 The calculation of the pseudoinverse $\inverse{(\bfC \bfZ)}$ can be reduced to solving a linear  least-squares problem and therefore its computing can be performed with the help of either the QR factorization or the SVD of the matrix $\bfC \bfZ$.

Although Algorithm~\ref{alg:proj_calc} for calculating the projection can be applied to the case of positive semidefinite weight matrices, the general case of the Cholesky factorization of a degenerate matrix $\bfW$ is complicated, see \cite[p. 201]{higham2002accuracy}.

\begin{algorithm}
	\caption{Calculation of $\winverse{\bfZ}{\bfW}$ and $\Proj_{\bfZ, \bfW} \bfx$ with the use of $\Sigminus = \bfC^\rmT \bfC$}
	\label{alg:proj_calc}
    \Input{$\bfZ \in \spC^{N\times r}$, $\bfW \in \spR^{N\times N}$ and $\bfx \in \spC^N$.}
	\begin{algorithmic}[1]
		\State{Compute the vector $\bfC \bfx$ and the matrix $\bfC \bfZ$.}
		\State{Calculate $\bfq = \inverse{(\bfC \bfZ)} (\bfC \bfx)$.}
		\State\Return{$\winverse{\bfZ}{\bfW} = \bfq \in \spR^{r\times N}$ and $\Proj_{\bfZ, \bfW} \bfx = \bfZ \bfq  \in \spR^{N}$.}
	\end{algorithmic}
\end{algorithm}

	However, there is a particular case of degenerate weight matrices, which correspond to a time series with missing values. In this case, the weight matrix $\bfW$ has zero columns and rows corresponding to missing entries and thereby can be easily processed. Denote $\bfu  \in \spR^N$ the vector with ones at the places of observations and zeros at the places of missing values, $\bfU = \diag(\bfu)$. Then the matrix $\bfW$ can be expressed as $\bfW = \bfU^\rmT \bfW_0 \bfU$. Suppose that $\bfW_0$ is positive definite. Consider the Cholesky decomposition $\bfW_0 = \bfC_0^\rmT \bfC_0$, where $\bfC_0$ is upper triangle, and set $\bfC = \bfC_0 \bfU$.
Then $\bfW = \bfC^\rmT \bfC$.
Note that if $\bfC_0$ is upper triangular with $p$ nonzero superdiagonals, then $\bfC$ is also upper triangular and has $p$ nonzero superdiagonals.

\section{Comparison of the VPGN and MGN algorithms}
\label{sec:comparison}
Let us compare the VPGN method and the proposed MGN method from the computational viewpoint.

\subsection{Design of comparison} \label{sec:allalgorithms}
We consider several versions of the algorithms VPGN and MGN with the directions of line search given by Algorithm~\ref{alg:gauss_newton_vp_step} and Algorithm~\ref{alg:gauss_newton_our_step} respectively.
The method VPGN is suggested in \cite{Usevich2014}.
The difference between VPGN and S-VPGN is described in \cite{Zvonarev.Golyandina2021}; in S-VPGN, the projections $\Proj_{\calZ(\bfa), \bfW}$ are calculated in a more stable way using the techniques introduced in \cite{Zvonarev.Golyandina2021} (see Algorithm~\ref{alg:fourier_basis_A} for calculating the special basis of $\calZ(\bfa)$ and Algorithm~\ref{alg:proj_calc} for calculating the projection itself in the case of a banded $\bfW$).
In MGN, the projections are always stably calculated. The addition of `-H' at the end of the method abbreviations means that the Compensated Horner scheme is used for calculating the projections.

The source code (\code{R} and \code{C++}) for these algorithms can be found in
\cite{Zvonarev2019}, where
the MGN and VPGN methods are implemented for the case of a common (not necessarily diagonal) weight matrix $\bfW$. The implementation of VPGN in
\cite{Zvonarev2019} extends that of \cite{Markovsky.Usevich2014}, which is suitable for diagonal weight matrices only, and has the same order of computational cost.

The comparison of VPGN, S-VPGN and S-VPGN-H was performed in \cite{Zvonarev.Golyandina2021}. This comparison demonstrated that the technique for the projection calculation, which was proposed in \cite{Zvonarev.Golyandina2021}, improves the algorithm numerical stability; however, this techniques cannot be applied to the calculation of $\Delta_k$ given in \ref{eq:delta_VPGN} and therefore does not allow one to considerably improve the VPGN method. The MGN method changes the iteration step in such a way (see \ref{eq:delta_MGN}) that the techniques developed in \cite{Zvonarev.Golyandina2021} can be used for all steps of the optimization algorithm.

\subsection{Theoretical comparison}
We start the comparison from comparing the algorithms by the computational costs. Then, we will compare the stability of the algorithms when the algorithms are comparable by the computational costs.
This depends on the structure of the weight matrix $\bfW$.
The special case of interest is the case when the weight matrix  $\Sigminus$ is $(2p+1)$-diagonal with a small $p$ (this is the case of autoregressive noise of order $p$ and therefore a natural assumption). Note that the special case when both $\Sigminus$
and $\Sigminus^{-1}$ are banded corresponds to the case of a diagonal matrix $\Sigminus$.

\subsubsection{Computational cost}
\label{sec:comp_cost}
Let us estimate computational costs in flops and study the asymptotic costs of one iteration as $N \rightarrow \infty$.

\paragraph{The MGN method}

The computation cost of the MGN method differs from that of the S-VPGN method only by the cost of calculating $\bfF_{\bfa^{(k)}}$ instead of $\bfJ_{S_\tau^\star}(\Ai0^{(k)})$ (compare Algorithms~\ref{alg:gauss_newton_our_step} and~\ref{alg:gauss_newton_vp_step}), since the projection $\Proj_{\calZ(\bfa^{(k)}), \bfW}$ is calculated by the same algorithms (with the cost $O(r N \log N + N r^2+Npr)$ \cite{Zvonarev.Golyandina2021}).

As it is shown in \cite{Zvonarev.Golyandina2021}, if $\Sigminus^{-1}$ is $(2p+1)$-diagonal, the asymptotic cost of the calculation of $\bfJ_{S_\tau^\star}(\Ai0^{(k)})$ is $O(N r^2 + N p^2)$ flops;
for the case when $\Sigminus$ is $(2p+1)$-diagonal, where $p > 0$, it takes $O(N^3)$.
The calculation of $\bfF_{\bfa^{(k)}}$ is performed by Algorithm~\ref{alg:fourier_grad}, which calls Algorithm~\ref{alg:fourier_basis_A}. The S-VPGN method contains Algorithm~\ref{alg:fourier_basis_A} with the cost $O(r N \log N + N r^2)$ for any form of $\bfW$ \cite{Zvonarev.Golyandina2021}; therefore, we need to discuss the cost of additional steps of Algorithm~\ref{alg:fourier_grad} only.

  Algorithm~\ref{alg:fourier_grad} includes the first two steps of Algorithm~\ref{alg:fourier_basis_A};
  the other steps can be performed in $O(r N \log N)$ flops.
  Thus, the cost of Algorithm~\ref{alg:fourier_grad} is $O(r N \log N + N r^2)$.
  Therefore, the asymptotic cost of one iteration of MGN given in Algorithm~\ref{alg:gauss_newton_our_step} is $O(N r^2 + Np^2 + r N \log N)$, or $O(Np^2 + N \log N)$ for a fixed rank $r$.

Table~\ref{table:costs} summarizes the computational costs for different methods.

\begin{table}[!hbt]
\begin{center}
	\caption{$r$-Rank approximation: Asymptotic computational costs in flops}
    \label{table:costs}
\begin{tabular}{|l|l|l|}
\hline
Method&$\bfW$ is $(2p+1)$-diagonal, $p>0$&$\bfW^{-1}$ is $(2p+1)$-diagonal\\
\hline
VPGN& $O(N^3)$ & $O(N r^2 + Np^2)$\\
S-VPGN,
S-VPGN-H&$O(N^3)$ & $O(N r^2 + Np^2 + r N \log N)$\\
MGN,
MGN-H& $O(N r^2 + Np^2 + r N \log N)$& $O(N r^2 + Np^2 +r N \log N)$\\
\hline
\end{tabular}
\end{center}
\end{table}

\begin{remark}
\label{rem:compcost}
Thus, if the inverse $\Sigminus^{-1}$ of the weight matrix $\Sigminus$ is $(2p+1)$-diagonal, then the computational cost of the proposed MGN method is slightly larger in comparison with the VPGN method.
However, if the weight matrix  $\Sigminus$ is $(2p+1)$-diagonal and $p>0$ (the case of autoregressive noise of order $p$),
then the computational cost of the MGN method is significantly smaller by order. In the case of a diagonal matrix $\bfW$ and a fixed $r$, the costs of MGN and VPGN are $O(N \log N)$ and $O(N)$ respectively.
\end{remark}

\subsubsection{Stability}
\label{sec:comp_thstab}

It is discussed in \cite{Zvonarev.Golyandina2021} that the main ``stability bottlenecks'' of the VPGN algorithms consists of solving the systems of linear equations. In fact, we say about inverting the matrices depending on $\bfa$, that is, on the coefficients of GLRR($\bfa$) governing the signal. The stability is measured by the  orders of the condition numbers of these matrices as the time-series length $N$ tends to infinity.

 For the MGN algorithm, the matrix $\bfA_g$ is inverted in Algorithm~\ref{alg:fourier_basis_A}. It was shown in \cite{Zvonarev.Golyandina2021} that the order of the condition number of the matrix $\bfA_g$ is $N^{t}$, where $t$ is the maximal multiplicity of the roots of the characteristic polynomial  $g(\bfa)$ \eqref{eq:pol_z} on the unit circle. It is worth to mention that the use of the compensated Horner scheme increases the accuracy of computing the diagonal elements of $\bfA_g$ and does not change its condition number.

For stability comparison, we consider the case when $\bfW^{-1}$ is banded, since otherwise the computational cost of the VPGN algorithm is very large. Independent noise with a diagonal covariance matrix $\bm\Sigma$ corresponds to this case in practice; then $\bfW = \bm\Sigma^{-1}$ is diagonal.
It was discussed in \cite{Zvonarev.Golyandina2021} with reference to \cite[Section 6.2]{Usevich2014} that the condition number of the ``bottleneck'' matrix of the VPGN algorithm has order $N^{2t}$.

Therefore, the advantage of the MGN algorithm in the accuracy is more significant in the case of the existence of unit-modulus roots of large enough multiplicity. In the general form \eqref{eq:model}, the multiplicity of a root is determined by the order of the corresponding polynomial $P_{m_k}$ and is equal to $m_k +1$. The root has unit modulus if $\alpha_k=0$. Therefore, the typical examples are polynomial signals and sinusoids with polynomial amplitude modulations.


\subsection{Numerical comparison} \label{subsec:speed}
For the numerical comparison of the considered algorithms, we use the example of a quadratic signal constructed in \cite{Zvonarev.Golyandina2021}, where a local solution of \eqref{eq:wls} is known.
Let us briefly describe this example.  For constructing a solution of rank $r=3$, the well-known theory explaining the relation of linear recurrence relations,
characteristic polynomials, their roots and the explicit form of the series, see e.g. the book \cite[Section 3.2]{Golyandina2013} with a brief description of this relation in the context of time series structure.

Let $\tsY_N^\star = (b t_1^2, \ldots, b t_N^2)^\rmT$, where $t_i$, $i=1,\ldots,N$, form the equidistant grid in $[-1; 1]$ and the constant $b$ is such that $\|\tsY_N^\star\|=1$. The series $\tsY_N^\star$ satisfies the GLRR($\bfa^*$) for $\bfa^* = (1, -3, 3, -1)^\rmT$. Since the last component of $\bfa^*$ is equal to $-1$, we can say that the series satisfies the LRR($\bfa^*$); the characteristic polynomial of this LRR has one unit root of multiplicity $t=3$.
Denote $\widehat \tsR_N = (c |t_1|, \ldots, c |t_N|)^\rmT$, where the constant $c$ is such that $\|\widehat \tsR_N\| = 1$.
Construct the observed series as $\tsX_N = \tsY_N^\star + \tsR_N$, where $\tsR_N = \widehat \tsR_N - \Proj_{\calZ((\bfa^*)^2), \Sigminus} \widehat \tsR_N$.
Thus, the pair $\tsX_0 = \tsY_N^\star$ and $\tsX = \tsX_N$ satisfies the necessary conditions for local minima \cite{Zvonarev.Golyandina2021}. The sufficient condition was tested numerically for $N < 100$.
Details of the example implementation see in \cite{Zvonarev.Golyandina2021}.

All computations are performed in double precision. When the compensated Horner scheme is used, the accuracy of calculations in double precision is the same as if the calculation would be in quadruple precision without the compensated Horner scheme.

\subsubsection{Stability} \label{subsec:basisacc}

The comparison is performed for the methods VPGN, S-VPGN-H, MGN and MGN-H for different $N$ from $100$ to $50000$.
For simplicity, consider the non-weighted case, when $\Sigminus$ is the identity matrix.

\begin{figure}[!hbt]
	\centering
	\includegraphics[width=10cm]{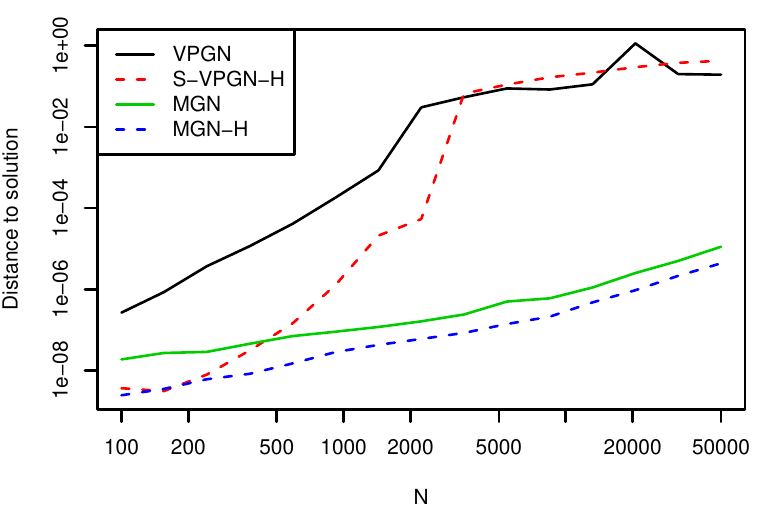}
\caption{Comparison of algorithms by distance to the solution, for different $N$.}
	\label{fig:VPGN_comp_disp}
\end{figure}

\begin{figure}[!hbt]
	\centering
	(a)\includegraphics[scale=1.00]{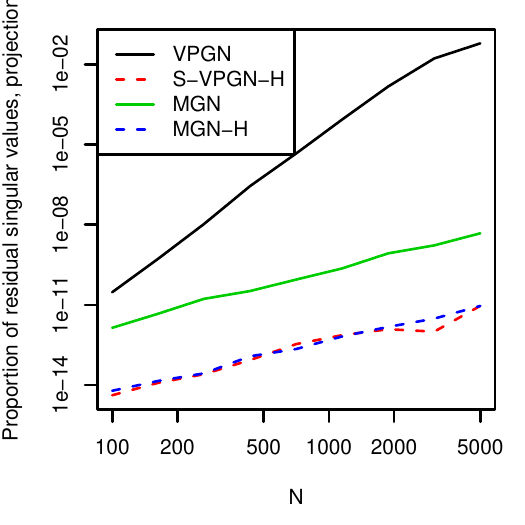}
	(b)\includegraphics[scale=1.00]{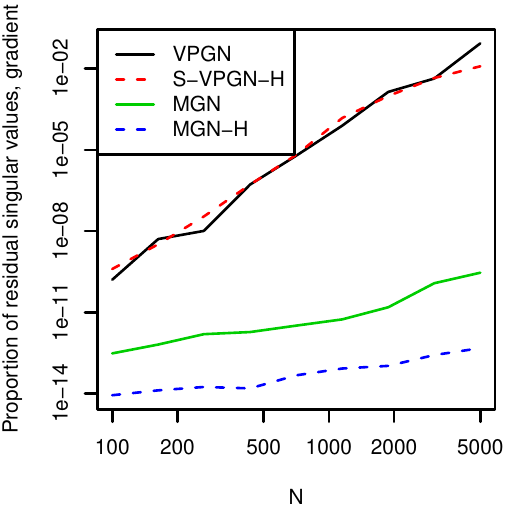}
	\caption{Shares of the residual eigenvalues, for different $N$; (a) the projection to $\overline{D_r}$ and (b) the Jacobian matrix with columns from a $2r$-dimensional space.}
	\label{fig:projection}
\end{figure}

 Denote $\widetilde \tsY^\star$ the result of an algorithm participating in the comparison.
We compare the algorithms by the accuracy, that is, by the Euclidean distance between $\widetilde \tsY^\star$ and $\tsY^\star_N$ (Fig.~\ref{fig:VPGN_comp_disp}).

The algorithms were started from the GLRR($\bfa_0$), where  $\bfa_0 =\bfa^* + 10^{-6} \bfd^\rmT$ and
each components of $\bfd$ is randomly distributed in $[-1,1]$.
We used 100 simulations to obtain the averaged results.
The iterations stop when the backtracking line search method (in step 6 of Algorithm~\ref{alg:common_scheme}) with decreasing search steps $\gamma = 1, 1/2, 1/4, \ldots, 2^{-50}$  does not decrease the value of the objective function.

Figure~\ref{fig:VPGN_comp_disp} shows that the accuracy of MGN and MGN-H is asymptotically much better than the accuracy of VPGN and S-VPGN-H; moreover, the errors of VPGN and S-VPGN-H increases drastically for large $N$. Note that in exact arithmetic, VPGN and S-VPGN-H would produce the same results; the same is true for the pair of MGN and MGN-H.

Figure~\ref{fig:projection} explains the difference between the algorithms. In \cite{Zvonarev.Golyandina2021}, the projection algorithms for calculating $\Proj_{\bfZ, \bfW} \tsX$ are compared by the difference with the true value. In Figure~\ref{fig:projection}(a), the comparison of the computed projectors to the set of time series of rank not larger than $r$ is performed by the proximity of their results to $r$-rank time series. Let $\tsY_N$ be the numerical result of the projection, choose $L=[N/2]$ and denote $\{\sigma_i\}_{i=1}^L$ the singular values of the $L$-trajectory matrix $\bfY=\calT_L(\tsY_N)$. We measure the projector accuracy as $\sqrt{\sum_{i=r+1}^L \sigma_i^2}/ \|\bfY\|_\rmF$, where $\|\bfY\|_\rmF^2 = \sum_{i=1}^L \sigma_i^2$ is the Frobenius norm. The projection $\Proj_{\bfZ, \bfW} \tsX$ is performed in the same way in the algorithms S-VPGN-H and MGN-H. Therefore, their accuracies coincide. MGN does not use the Horner compensation scheme and therefore it is a bit worse. The algorithm used in VPGN is less numerically accurate.

The key difference between the VPGN and MGN algorithms is in the 5-th step of Algorithm~\ref{alg:common_scheme}, that is, in the way of calculation of the Jacobian matrix; this can be performed in MGN with improved accuracy. Since each column of the Jacobian matrix (that is, each gradient vector) has rank $2r$ (Theorem~\ref{th:tangent}), we can check the accuracy of its calculation by summing the shares $\sqrt{\sum_{i=2r+1}^L \sigma_i^2}/ \|\bfY\|_\rmF$ for columns $\tsY_N$ of $\bfJ_{S_{\tau}^\star}(\Ai0^{(k)})$ for VPGN and of $\big(\bfI_N - \Proj_{\calZ(\fullop(\Ai0^{(k)})), \bfW}\big)\widehat \bfF_{\bfa}$ for MGN, see \eqref{eq:delta_VPGN} and \eqref{eq:delta_MGN} and  Algorithms~\ref{alg:gauss_newton_vp_step} and \ref{alg:gauss_newton_our_step}.

Figure~\ref{fig:projection}(b) confirms that VPGN and S-VPGN-H lack numerical stability when computing the direction (step 5 of Algorithm~\ref{alg:common_scheme}); recall that the multiplicity $t$ equals 3 for this example (see discussion in Section~\ref{sec:comp_thstab}).
 Figure~\ref{fig:projection}(a) shows that VPGN produces a resultant series far from $r$-rank one, for large $N$. Both of these factors lead to that VPGN and S-VPGN-H perform worse than the MGN and MGN-H methods (a ``plateau'' visible in Figure~\ref{fig:VPGN_comp_disp} for VPGN and S-VPGN-H for $N > 2000$ is due to the fact that the distance can not be greater than the series norm).

\subsubsection{Computational cost}
For effectively implemented algorithms, the computational speed should have the same order as the theoretical computational cost in flops. Let us numerically confirm Remark~\ref{rem:compcost}.
We will consider the computational time for different implementations of Algorithms \ref{alg:gauss_newton_vp_step} and \ref{alg:gauss_newton_our_step}, which calculate the search directions. This time characterizes the computational speed
of one iteration.
We consider different time series lengths $N$ and two types of the weight matrix $\bfW$, the identity matrix and a 3-diagonal matrix, which is the inverse of the covariance matrix of AR(1) with the coefficient 0.9  (the explicit form of this inverse can be found e.g. in \cite{Golyandina2020a}). The speed is estimated with the help of the example described in Section~\ref{subsec:basisacc}.

\begin{figure}[!hbt]
	\centering
	(a)\includegraphics[scale=1.00]{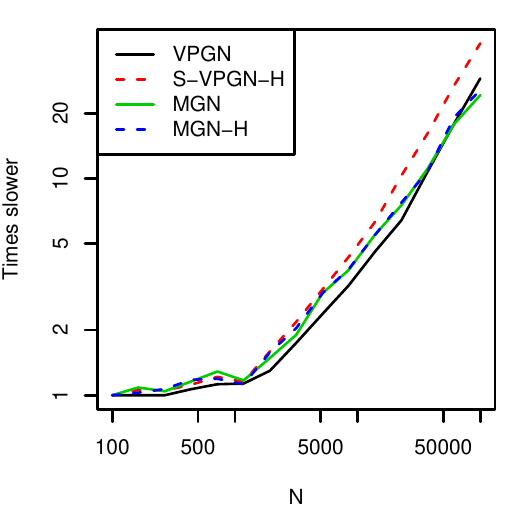}
	(b)\includegraphics[scale=1.00]{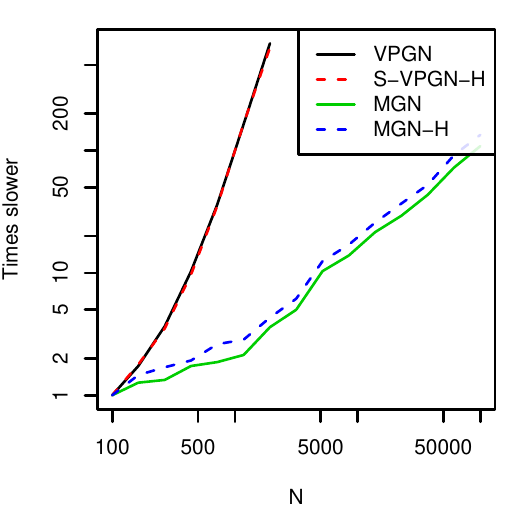}
	\caption{Comparison of algorithms by CPU times of one iteration for different $N$; (a) diagonal $\bfW$ and (b) 3-diagonal $\bfW$.}
	\label{fig:comp_time}
\end{figure}

The results for the CPU time are depicted in Fig.~\ref{fig:comp_time}. Since we compare asymptotic behavior (as $N\to\infty$), we eliminate the constant time, which does not depend on $N$, in the following way.
 For each algorithm, we consider the CPU times for different values of $N$ starting from 100 and then divide them by the CPU time for $N$ equal 100.
  Note that if $\bfW$ is diagonal, the computational times of the algorithms are asymptotically almost the same. However, if
  $\bfW$ contains three diagonals, the computational times for the methods MGN and MGN-H are much smaller than that for the methods VPGN and S-VPGN-H.

\subsection{Signal estimation using MGN: with and without gaps} \label{subsec:model}
Consider a time series $\tsY_{50}$ similar to the one considered in \cite{Ishteva.etal2014}, which is the sum of a signal of rank $r=4$ and Gaussian white noise. That is, let the signal $\tsS_{50}$ have the following form: $\tsS_{50} = (s_1, \ldots, s_{50})^\rmT$, where
\begin{equation}
\label{eq:ex_mv}
s_i = 0.9^i \cos\left(\frac{\pi}{5} i\right) + \frac{1}{5} 1.05^i \cos \left(\frac{\pi}{12} i + \frac{\pi}{4}\right), \quad i = 1, \ldots, 50,
\end{equation}
and
\begin{equation*}
\tsY_{50} = \tsS_{50} + 0.2 \frac{\tsR_{50}}{\|\tsR_{50}\|}\|\tsS_{50}\|;
\end{equation*}
here the series $\tsR_{50}$ consists of i.i.d. normal random variables with zero mean and unit variance. Note that since $\tsY_{50}$ contains a random component, we were not able to reproduce the time series studied in \cite{Ishteva.etal2014} exactly; we used \code{set.seed(15)} in \code{R version 4.0.2}, see \cite{Zvonarev2019}.

Let us consider two versions of the time series $\tsY_{50}$, the first one is without missing data and the second time series with artificial gaps at positions $i = 10 \ldots 19$ and $i = 35 \ldots 39$, and construct two estimates of
the signal by the MGN-H method.
Since the noise is white,
the identity matrix $\bfW = \bfI_{50}$ was taken for the case without gaps; for the case with gaps, we consider $\bfW_0 = \bfI_{50}$ and construct $\bfW$ changing ones on the diagonal of $\bfW_0$ at the positions of missing data to zeros (see Section \ref{sec:missing}).
For constructing the initial GLRR, we replace the missing entries by the mean value of the time series and then set
the GLRR coefficients equal to $(r+1)$-th left singular vector  of the SVD of the $(r+1)$-trajectory matrix
$\calT_{r+1}(\tsY_{50})$.

\begin{figure}[!hbt]
	\centering
	(a)\includegraphics[scale=1.00]{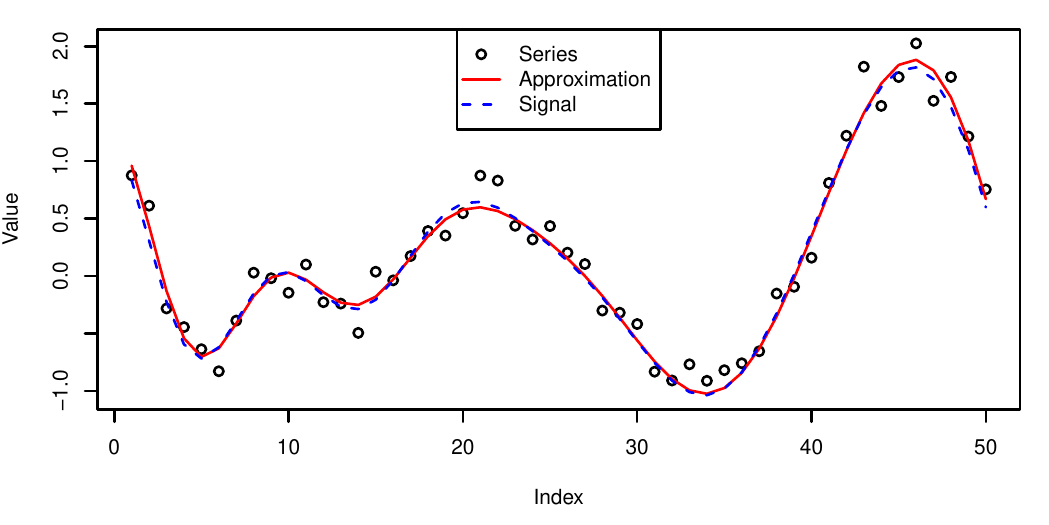}
	(b)\includegraphics[scale=1.00]{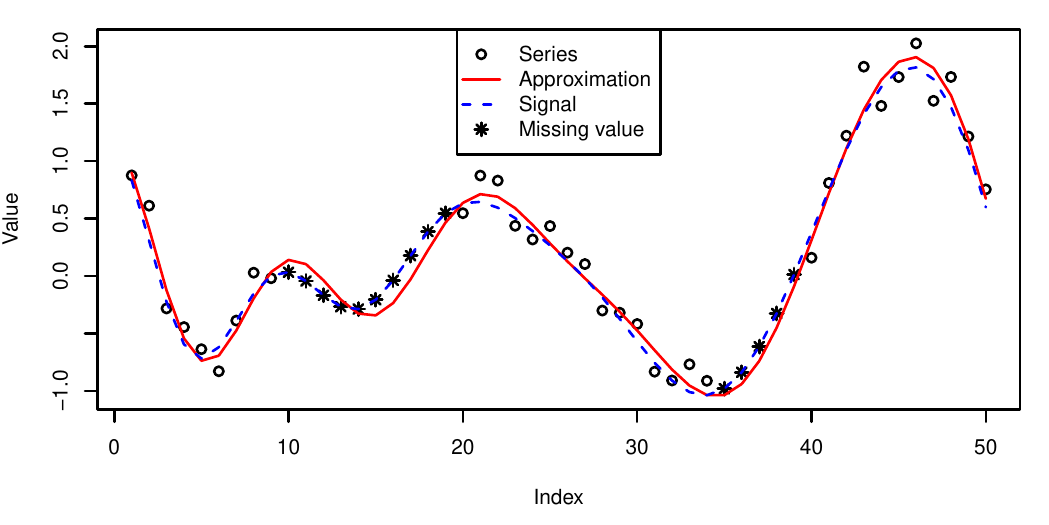}
	\caption{Estimates of the signal $\tsS_{50}$ using the MGN algorithm, (a) without gaps (b) with gaps.}
	\label{fig:model}
\end{figure}

The results are presented in Figure \ref{fig:model}. The series $\tsY_{50}$ is indicated by the black dots, the signal $\tsS_{50}$ is depicted by the blue line, and the obtained approximation $\widetilde \tsS$ is shown by the red solid line. Note that in both cases $\widetilde \tsS$ gives a fairly close estimate of $\tsS_{50}$, despite even a large gap at $10 \ldots 19$ in the case with missing values.

\medskip
Let us consider the same signal \eqref{eq:ex_mv} corrupted by AR(1) noise with coefficient 0.9 (red noise) instead of white noise. We compare two versions of the MGN-H algorithm, with $\bfW^{(\rmW)} = \bfI_{N}$ and with $\bfW^{(\rmR)} = \bm\Sigma^{-1}$, where $\bm\Sigma$ is the $N\times N$ covariance matrix of the AR process. The latter corresponds to the MLE estimate. To eliminate the influence of the initial GLRR on the comparison, we take it as the minimal GLRR governing the signal.

The estimates of RMSE are obtained by 1000 iterations; in the parentheses we indicate the computational time in seconds. For signal estimation without gaps, we have RMSE equal to 0.075 for  $\bfW=\bfW^{(\rmW)}$ (80 sec) and 0.066 for $\bfW=\bfW^{(\rmR)}$ (75 sec). For gap imputation accuracy in presence of gaps, we have 0.136 for $\bfW_0=\bfW^{(\rmW)}$ (105 sec) and 0.097 for $\bfW_0=\bfW^{(\rmR)}$ (100 sec). The difference between the errors is statistically significant. The computational times are similar and even a bit smaller for the MLE version. This confirms that MGN-H is well working for time series with possible gaps and autoregressive noise.

\section{Conclusion}
\label{sec:conclusion}
In this paper we presented a new iterative algorithm (MGN, Algorithm~\ref{alg:common_scheme} with $\Delta_k$ calculated by Algorithm~\ref{alg:gauss_newton_our_step}) for computing the numerical solution to the
problem~\eqref{eq:wls} and compared it with a state-of-art algorithm based on the variable projection approach (VPGN, Algorithm~\ref{alg:common_scheme} with $\Delta_k$ calculated by Algorithm~\ref{alg:gauss_newton_vp_step}).
We showed that the proposed algorithm MGN allows the implementation, which is more numerically stable especially for the case of
multiple unit-modulus roots of the characteristic polynomial (in particular, for polynomial series, where the multiplicity
is equal to the polynomial order plus one). This effect can be explained by the inversion of matrices with
condition number of order $N^{t}$ in MGN, where $t$ is the multiplicity, while the direct implementation
of VPGN deals with matrices with condition number of order $N^{2t}$.

The comparison of computational costs in Section~\ref{sec:comp_cost} shows that the algorithm MGN has slightly larger costs for the case of banded inverses $\bfW^{-1}$ of weight matrices. However, in the case of autoregressive noise with covariance matrix $\bm\Sigma$,
the corresponding weight matrix $\bfW = \bm\Sigma^{-1}$, which provides the MLE estimate, is banded itself and $\bfW^{-1}$ is not banded.
Then the proposed algorithm MGN has a much lower computational cost in comparison with VPGN.

An important feature of the MGN algorithm is that it can be naturally extended to the case of missing data
without increasing the computational cost (see Sections~\ref{sec:missing} and \ref{subsec:model}).


\section*{Acknowledgments} The reported study was funded by RFBR, project number 20-01-00067.

\bibliography{references}

\section*{Author Biography}

\begin{biography}{\includegraphics[width=66pt,height=74pt]{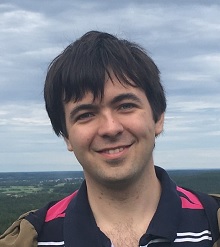}}{\textbf{Nikita Zvonarev.} Nikita Zvonarev is an Assistant Professor at St. Petersburg State University, St.
Petersburg, Russia. Dr. Zvonarev received the M.S. degree in Applied Mathematics in 2015 and Ph.D. degree in Computational Mathematics with the thesis entitled ``Structured time series approximations'' in 2018 from St. Petersburg State University. His research interests are optimization methods and time series analysis.}
\phantom{q q q q q q q q q q q q q q q q q q q q q q q q q q q q q q q q q q q q q q q q q q q q q q q q q q q q q q q q q q q q q q q q q q q q q q q q q q q q q q q q q q q q q q q q q q q q q q q q q}
\end{biography}

\begin{biography}{\includegraphics[width=66pt,height=86pt]{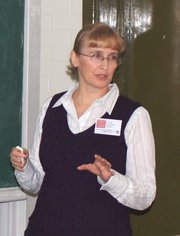}}{\textbf{Nina Golyandina.} Nina Golyandina is an Associate Professor at St. Petersburg State University, St.
Petersburg, Russia. Dr. Golyandina received the M.S. degree in Applied Mathematics in 1985 and Ph.D. degree in Probability Theory and Mathematical Statistics
in 1998 from St. Petersburg State University. Her research interests are in applied
statistics and time series analysis. She is the co-author of three monographs devoted to singular spectrum analysis.}
\end{biography}

\end{document}